\documentclass[11pt, a4paper]{article}

\usepackage[T1]{fontenc}
\usepackage[ansinew]{inputenc}
\usepackage{amsmath}
\usepackage{amsfonts}
\usepackage{txfonts}
\usepackage[amsmath, standard, thmmarks, thref]{ntheorem}
\usepackage{fancyhdr}
\usepackage{textcomp}
\usepackage{stmaryrd}
\usepackage{graphicx}
\usepackage{psfrag}
\usepackage{float}

\newtheorem{Thm}{Theorem}[section]
\newtheorem{Le}[Thm]{Lemma}

\newtheorem{Rem}[Thm]{Remark}

\theoremstyle{plain}
\theoremheaderfont{\normalfont\bfseries}
\theorembodyfont{}

\theoremstyle{nonumberplain}
\theoremheaderfont{\itshape}
\theorembodyfont{}
\theoremsymbol{$\square$}
\newtheorem{Pro}{Proof:}

\DeclareMathOperator{\diam}{diam}

\DeclareMathOperator{\graph}{graph}
\DeclareMathOperator{\C}{C}
\DeclareMathOperator{\D}{D}

\DeclareMathOperator{\dist}{dist}

\DeclareMathOperator{\spt}{spt}

\DeclareMathOperator{\genus}{genus}

\DeclareMathOperator{\Rdrei}{\varmathbb{R}^3}
\DeclareMathOperator{\Rzwei}{\varmathbb{R}^2}
\DeclareMathOperator{\Reins}{\varmathbb{R}}
\DeclareMathOperator{\Leins}{\mathcal{L}^1}
\DeclareMathOperator{\Lzwei}{\mathcal{L}^2}
\DeclareMathOperator{\Ldrei}{\mathcal{L}^3}
\DeclareMathOperator{\Heins}{\mathcal{H}^1}
\DeclareMathOperator{\Hzwei}{\mathcal{H}^2}

\DeclareMathOperator{\Will}{\mathcal{W}}
\DeclareMathOperator{\cur}{\mathcal{M}}
\DeclareMathOperator{\M}{\mathcal{M}\,_\sigma}

\DeclareMathOperator{\schi}{\chiup_{_k}}
\DeclareMathOperator{\B}{\mathcal{B}_\varepsilon}
\DeclareMathOperator{\diver}{div}
\DeclareMathOperator{\A}{\textnormal{A}}
\DeclareMathOperator{\Be}{\textnormal{B}}
\DeclareMathOperator{\Ce}{\textnormal{C}}

\begin{document}

\parindent+0cm
\numberwithin{equation}{section}
\renewcommand{\thefootnote}{}

\begin{center}
\Huge{Willmore minimizers with prescribed isoperimetric ratio}
\end{center}
\vspace{0,5cm}
\begin{center}
Johannes Schygulla*\\Mathematisches Institut der Albert-Ludwigs-Universit\"at Freiburg\\Eckerstra\ss e 1, D-79104 Freiburg, Germany\\email: johannes.schygulla@math.uni-freiburg.de
\end{center}

\footnotetext{*J.Schygulla was supported by the DFG Collaborative Research Center SFB/Transregio 71.}

\begin{center}

\quote{{\bf Abstract:} Motivated by a simple model for elastic cell membranes, we minimize the Willmore functional among two-dimensional spheres embedded in $\Rdrei$ with prescribed isoperimetric ratio.
\bigskip

{\bf Key words:} Willmore functional, geometric measure theory. 

\bigskip

{\bf MSC:} 53 A 05, 49 Q 20, 49 Q 15, 74 G 65}

\bigskip

\end{center}

\section{Introduction}
In the spontaneous curvature model for lipid bilayers due to Helfrich \cite{H}, the membrane of a vesicle is described as a two-dimensional, embedded surface $\Sigma\subset\Rdrei$, whose energy is given by
\begin{equation*}
E(\Sigma)=\kappa\int_\Sigma(H-C_0)^2\,d\mu+\kappa_G\int_\Sigma K\,d\mu,
\end{equation*}
where $H$, $K$ denote the mean curvature and Gauss curvature, $\mu$ is the induced area measure and $\kappa$ and $\kappa_G$ are constant bending coefficients.\\[0,3cm]
Restricting to surfaces $\Sigma$ of the type of the sphere, the second term reduces to the constant $4\pi\kappa_G$ by the Gauss-Bonnet theorem. Reducing further to the simplest case of spontaneous curvature $C_0=0$, the energy becomes up to a factor the Willmore energy
\begin{equation}\label{I2}
\Will(\Sigma)=\frac{1}{4}\int_\Sigma|\vec{H}|^2\,d\mu.
\end{equation}
According to \cite{H}, the shapes of the vesicles should be minimizers of the elastic energy $E$ subject to prescribed area and enclosed volume. Since the Willmore energy is scaling invariant, the two constraints actually reduce to the condition that the isoperimetric ratio of the surface $\Sigma$, given by
\begin{equation}\label{I3}
I(\Sigma)=\left(6\sqrt{\pi}\,\right)^{\frac{1}{3}}\frac{V(\Sigma)^\frac{1}{3}}{A(\Sigma)^\frac{1}{2}},
\end{equation}
is prescribed. Here $A(\Sigma)$ denotes the area of $\Sigma$ and $V(\Sigma)$ the volume enclosed by $\Sigma$, i.e. the volume of the bounded component of $\Rdrei\setminus\Sigma$. The normalizing constant $(6\sqrt{\pi}\,)^{\frac{1}{3}}$ is chosen such that $I(\Sigma)\in(0,1]$, in particular $I(\varmathbb{S}^2)=1$. \\[0,3cm]
For given $\sigma\in(0,1]$, we denote by $\M$ the class of smoothly embedded surfaces $\Sigma\subset\Rdrei$ with the type of ${\varmathbb S}^2$ and with $I(\Sigma)=\sigma$, and we introduce the function
$$\beta:(0,1]\to\Reins_+,\quad\beta(\sigma)=\inf_{\Sigma\in\M}\Will(\Sigma).$$
We have $\mathcal{M}_1=\Big\{\text{round spheres}\subset\Rdrei\Big\}$ and $\beta(1)=4\pi$. \\[0,3cm]
Here we prove the following result.

\begin{Thm}\label{MainThm}
For every $\sigma\in(0,1)$ there exists a surface $\Sigma\in\M$ such that
$$\Will(\Sigma)=\beta(\sigma).$$
Moreover the function $\beta$ is continuous, strictly decreasing and satisfies
$$\lim_{\sigma\searrow0}\beta(\sigma)=8\pi.$$
\end{Thm}

Assuming axial symmetry, several authors computed possible candidates for minimizers by solving numerically the Euler-Lagrange equations (see \cite{BLS}, \cite{DH}). In \cite{NT} the authors prove existence of a one-parameter family of critical points bifurcating from the sphere. It appears that so far no global existence results for the Helfrich model have been obtained. In order to prove \thref{MainThm} we adopt the methods of L. Simon in \cite{SL}, where he proved existence of Willmore minimizers for fixed $\genus p=1$.\\[0,3cm]
Moreover we show the following result.

\begin{Thm}\label{Thm2}
Let $\{\sigma_k\}_{k\in\varmathbb{N}}\subset(0,1)$ such that $\sigma_k\to0$ and $\Sigma_k\in\mathcal{M}_{\sigma_k}$ such that $\Will(\Sigma_k)=\beta(\sigma_k)$. After translation and scaling (such that $0\in\Sigma_k$ and $\Hzwei(\Sigma_k)=1$), there exists a subsequence $\Sigma_{k'}$ which converges to a double sphere in the sense of measures, namely
$$\mu_{k'}\to\mu\quad\text{in }C_c^0(\Rdrei)',$$
where $\mu_{k'}=\Hzwei\llcorner\Sigma_{k'}$ and $\mu=2\Hzwei\llcorner\partial B_r(a)$ for some $r>0$ and $a\in\Rdrei$.
\end{Thm}

We now briefly outline the content of the paper. In section 2 we prove that $\beta$ is decreasing and $\beta(\sigma)<8\pi$ for all $\sigma\in(0,1]$. In section 3 we prove \thref{MainThm}. Section 4 is dedicated to the proof of \thref{Thm2} using similar techniques as in the proof of \thref{MainThm}. Finally in the appendix we collect some important results we need during the proofs, as for example the graphical decomposition lemma and the Monotonicity formula proved by Simon in \cite{SL}.\\[0,3cm]
This work was done within the framework of project B.3 of the DFG Collaborative Research Center SFB/Transregio 71. I would like to thank my advisor Prof. Ernst Kuwert for his support. I also would like to express my gratitude for the support I received from the DFG Collaborative Research Center SFB/Transregio 71.

\section{Upper bound for the Infimum}
In this section we prove an upper bound for the infimum of the Willmore energy in the class $\M$. The proof is based on the inversion of a catenoid at a sphere together with an argument involving the Willmore flow and its properties. A reference where the authors also analyze inverted catenoids and their relation to the Willmore energy is \cite{CVG}. For the part concerning the Willmore flow see \cite{KS}.
\begin{Le}\label{48}
The function $\beta$ is decreasing and
$$\beta(\sigma)=\inf_{\Sigma\in\M}\Will(\Sigma)<8\pi\quad\text{for all }\sigma\in(0,1].$$
\end{Le}
\textit{Proof: }Define the (scaled) catenoid in $\Rdrei$ as the image of $g_a:\Reins\times[0,2\pi)\to\Rdrei$ given by
\begin{equation*}
g_a(s,\theta)=\left(a\cosh\frac{s}{a}\cos\theta, a\cosh\frac{s}{a}\sin\theta, s\right),
\end{equation*}
where $a>0$ is a positive constant. Next we invert this catenoid at the sphere $\partial B_1(e_3)$ to get the function $f_a=I\circ g_a$, where $I(x)=e_3+\frac{x-e_3}{|x-e_3|^2}$ describes the inversion at the sphere. Define the set $\Sigma_a\subset\Rdrei$ by 
$$\Sigma_a=f_a\,\Big(\Reins\times[0,2\pi)\,\Big)\cup\{e_3\}.$$
\vspace{-0,5cm}
\begin{figure}[h]
\begin{center}
\psfrag{A}{\hspace{-0,35cm}\small{graph $u_+$}}
\psfrag{B}{\hspace{-0,4cm}\small{graph $u_-$}}
\psfrag{C}{\small{$\mathcal{U}$}}
\includegraphics[scale=0.5]{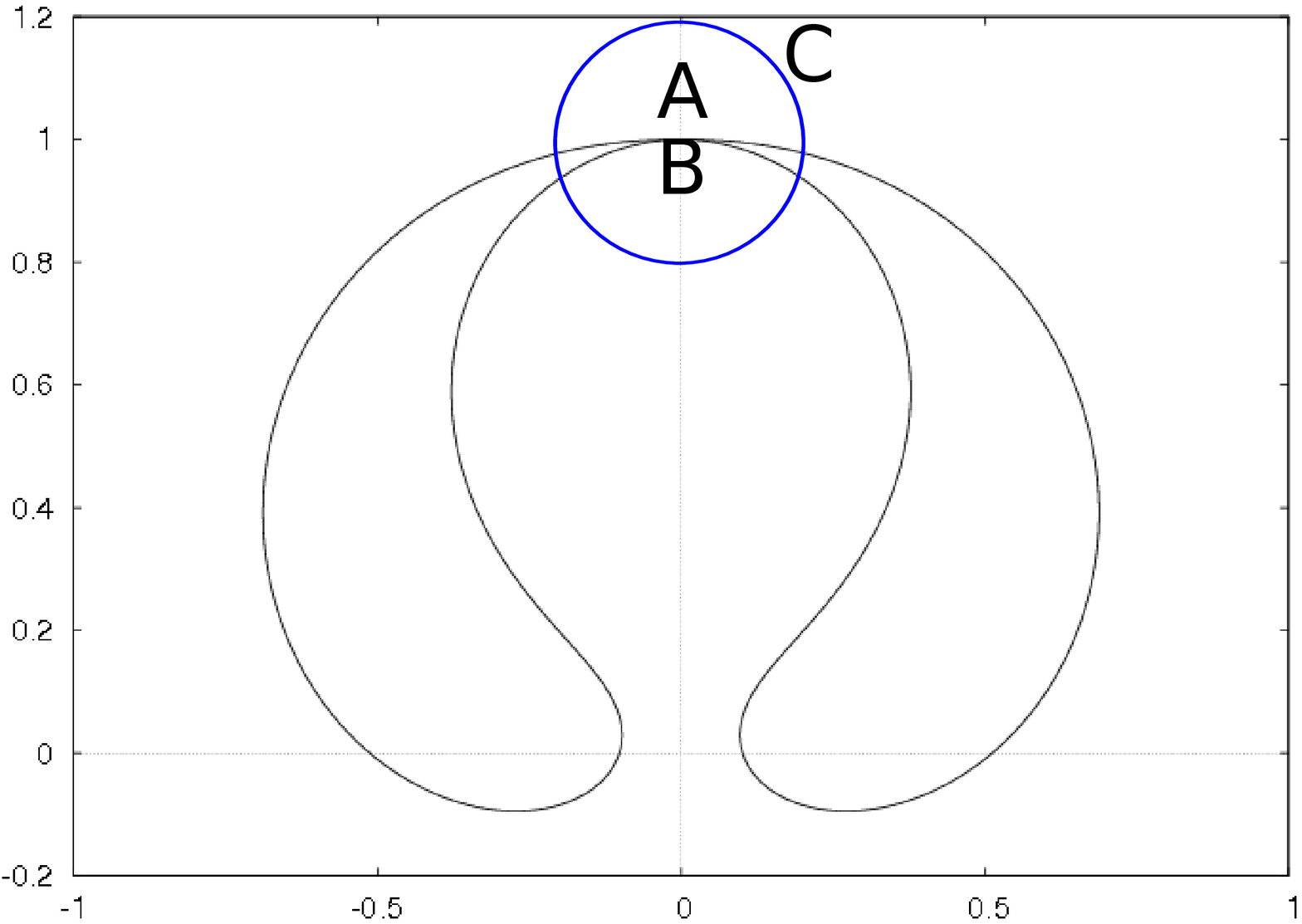}
\caption{$f_a((\Reins\times\{0\})\cup(\Reins\times\{\pi\}))$ for $a=0.1$. $\Sigma_a$ results from rotation.}
\end{center}
\end{figure}

First of all $\Sigma_a$ is smooth away from $e_3$. Because of the inverse function theorem and by explicit calculation there exists an open neighborhood $\mathcal{U}$ of $e_3$ in which $\Sigma_a$ can be written as $\graph u_+\cup\graph u_-$, where $u_{\pm}\in C^{1,\alpha}(B_R(0))\cap W^{2,p}(B_R(0))$ for all $\alpha\in(0,1)$, $p\ge1$ and some $R>0$, and which are smooth away from the origin. Moreover direct calculation yields $\Will(\Sigma_a)=8\pi$. Since variations of $\Sigma_a$ away from $e_3$ correspond to variations of the catenoid away from infinity and since the catenoid is a minimal surface, the $L^2$-gradient $\vec{\Will}(f_a)$ of the Willmore energy of $f_a$ satisfies
$$\vec{\Will}(f_a)=0\quad\text{on }(-\infty,\infty)\times[0,2\pi).$$
Since $u_{\pm}\in C^{1,\alpha}(B_R(0))\cap W^{2,p}(B_R(0))$ for all $\alpha\in(0,1)$ and $p\ge1$, it follows that the $L^2$-gradient of the Willmore energy of $\graph u_\pm$ satisfies
$$\vec{\Will}(F_\pm)=0\quad\text{on }B_R(0)\setminus\{0\},$$
where $F_\pm(x,y)=(x,y,u_\pm(x,y))$. Let $\phi\in C_c^\infty(B_R(0))$ and define the function $F_\pm^t(x,y)=(x,y,u_\pm(x,y)+t\phi(x,y))$. For $\Omega\subset B_R(0)$ denote by $\Will(F_\pm^t,\Omega)$ the Willmore energy of $F_\pm^t$ restricted to $\Omega$. Because of the given regularity of $u_\pm$ and since $\spt\phi\subset\subset B_R(0)$ it follows that
$$\frac{d}{dt}\Will(F_\pm^t)_{|_{t=0}}=\frac{d}{dt}\Will(F_\pm^t,B_R(0))_{|_{t=0}}=\lim_{\varepsilon\to0}\frac{d}{dt}\Will(F_\pm^t,B_R(0)\setminus B_\varepsilon(0))_{|_{t=0}}.$$
Since $\vec{\Will}(F_\pm)=0$ on $B_R(0)\setminus B_\varepsilon(0)$, it follows from the first variation formula for the Willmore energy that only a boundary term remains. Exploiting this boundary integral yields
$$\lim_{\varepsilon\to0}\frac{d}{dt}\Will(F_\pm^t,B_R(0)\setminus B_\varepsilon(0))_{|_{t=0}}=\mp c\phi(0),$$
where $c>0$ is a positive constant. This shows that the first variation of the Willmore energy of $\graph u_+$ is negative for variations in the direction $e_3$ and that the first variation of the Willmore energy of $\graph u_-$ is negative for variations in the direction $-e_3$. Now notice that the isoperimetric ratio $I(\Sigma_a)\to0$ as $a\to0$ and that $\Sigma_a$ can be parametrized over $\varmathbb{S}^2$. After approximation by smooth surfaces we have therefore shown that for every $\varepsilon>0$ there exists a smooth, embedded surface $\Sigma\subset\Rdrei$ of the type of $\varmathbb{S}^2$, with isoperimetric ratio $I(\Sigma)<\varepsilon$ and $\Will(\Sigma)<8\pi$. Using Theorem 5.2 in \cite{KS}, the Willmore flow $\Sigma_t$ with initial data $\Sigma$ exists smoothly for all times and converges to a round sphere such that $\Will(\Sigma_t)$ is decreasing in $t$. This shows $\beta(\sigma)<8\pi$. In order to prove the monotonicity let $\sigma_0\in(0,1)$ and $\varepsilon>0$ such that $\beta(\sigma_0)+\varepsilon<8\pi$. Let $\Sigma_0\in\mathcal{M}_{\sigma_0}$ such that $\Will(\Sigma_0)\le\beta(\sigma_0)+\varepsilon$. Again the Willmore flow $\Sigma_t$ with initial data $\Sigma_0$ exists smoothly for all times, converges to a round sphere and $\Will(\Sigma_t)$ is decreasing in $t$. Therefore for every $\sigma\in(\sigma_0,1]$ there exists a surface $\Sigma\in\M$ with $\Will(\Sigma)\le\Will(\Sigma_0)\le\beta(\sigma_0)+\varepsilon$ and the lemma follows by letting $\varepsilon\searrow0$.\hfill$\square$

\section{Proof of \thref{MainThm}}
For $\sigma\in(0,1)$ let $\big\{\Sigma_k\big\}_{k\in\varmathbb{N}}\subset\M$ be a minimizing sequence. Since the Willmore energy is invariant under translations and scalings and in view of \thref{48} we may assume that for some $\delta_0>0$
\begin{equation}\label{49}
\Hzwei(\Sigma_k)=1\,,\quad0\in\Sigma_k\,,\quad\Will(\Sigma_k)\le8\pi-\delta_0.
\end{equation}

Using Lemma 1.1 in \cite{SL} we get an uniformly diameter bound for $\Sigma_k$ and therefore
\begin{equation}\label{2}
\Sigma_k\subset B_R(0)\quad\text{for some }R<\infty.
\end{equation}

Define the integral, rectifiable 2-varifolds $\mu_k$ in $\Rdrei$ by
\begin{equation}\label{3}
\mu_k=\Hzwei\llcorner\Sigma_k.
\end{equation}

By a compactness result for varifolds (see \cite{SLGMT}), there exists an integral, rectifiable 2-varifold $\mu$ in $\Rdrei$ with density $\theta(\mu,\cdot)\ge1$ $\mu$-a.e. and weak mean curvature vector $\vec{H}\in L^2(\mu)$, such that (after passing to a subsequence) $\mu_k\to\mu$ in $C_c^0(\Rdrei)'$ and
\begin{equation}\label{4''}
\lim_{k\to\infty}\int\big<X,\vec{H}_k\big>\,d\mu_k=\int\big<X,\vec{H}\big>\,d\mu\quad\text{for all }X\in C^1_c(\Rdrei,\Rdrei),
\end{equation}
\begin{equation}\label{4'}
\frac{1}{4}\int_U|\vec{H}|^2\,d\mu\le\liminf_{k\to\infty}\frac{1}{4}\int_U|\vec{H}_k|^2\,d\mu_k\le8\pi-\delta_0\quad\text{for all open }U\subset\Rdrei.
\end{equation}
\thref{Monotonie} and (\ref{4'}) applied to $U=\Rdrei$ yield
\begin{equation}\label{4}
\theta^2(\mu,x)\in\left[1,2-\frac{\delta_0}{4\pi}\right]\quad\text{for all }x\in\spt\mu.
\end{equation}
Since $\mu$ is integral we also get that
\begin{equation}\label{4'''}
\theta^2(\mu,x)=1\quad\text{for }\mu\text{-a.e. }x\in\spt\mu.
\end{equation}

Our candidate for a minimizer is given by
\begin{equation}\label{5}
\Sigma=\spt\mu.
\end{equation}

Using \thref{Monotonie} we get (up to subsequences) as in \cite{SL}, page 310, that
\begin{equation}\label{6}
\Sigma_k\to\Sigma\quad\text{in the Hausdorff distance sense.}
\end{equation}
Therefore (\ref{2}) and the varifold convergence yield that $\Sigma\subset\overline{B_R(0)}$ and $\mu(\Rdrei)=1$.\\[0,3cm]
In order to prove regularity we would like to apply Simon's graphical decomposition lemma \thref{Decomposition} to $\Sigma_k$ simultaneously for infinitely many $k\in\varmathbb{N}$. But the most important assumption in the graphical decomposition lemma is that the $L^2$-norm of the second fundamental form is locally small, which we will need simultaneously for infinitely many $k\in\varmathbb{N}$. Therefore we define the so called bad points with respect to a given $\varepsilon>0$ in the following way: Define the radon measures $\alpha_k$ on $\Rdrei$ by
$$\alpha_k=\mu_k\llcorner|\rm{A}_k|^2.$$
From the Gauss-Bonnet formula and (\ref{49}) it follows that $\alpha_k(\Rdrei)\le24\pi$. By compactness there exists a radon measure $\alpha$ on $\Rdrei$ such that (after passing to a subsequence) $\alpha_k\to\alpha$ in $C_c^0(\Rdrei)'.$ It follows that $\spt\alpha\subset\Sigma$ and $\alpha(\Rdrei)\le24\pi$. Now we define the bad points with respect to $\varepsilon>0$ by
\begin{equation}\label{7}
\B=\left\{\xi\in\Sigma\,\big|\,\alpha(\{\xi\})>\varepsilon^2\right\}.
\end{equation}

Since $\alpha(\Rdrei)\le c$, there exist only finitely many bad points. Moreover for $\xi_0\in\Sigma\setminus\B$ there exists a $0<\rho_0=\rho_0(\xi_0,\varepsilon)\le1$ such that $\alpha(B_{\rho_0}(\xi_0))<2\varepsilon^2$, and since $\alpha_k\to\alpha$ weakly as measures we get
\begin{equation}\label{9}
\int_{\Sigma_k\cap B_{\rho_0}(\xi_0)}|\rm{A}_k|^2\,d\Hzwei\le2\varepsilon^2\quad\text{for }k\text{ sufficiently large}.
\end{equation}

Now fix $\xi_0\in\Sigma\setminus\B$ and let $\rho_0$ as in (\ref{9}). Let $\xi\in\Sigma\cap B_\frac{\rho_0}{2}(\xi_0)$. We want to apply Simon's graphical decomposition lemma to show that the surfaces $\Sigma_k$ can be written as a graph with small Lipschitz norm together with some "pimples" with small diameter in a neighborhood around the point $\xi$. This is done in exactly the same way Simon did in \cite{SL}. We just sketch this procedure. By (\ref{6}) there exists a sequence $\xi_k\in\Sigma_k$ such that $\xi_k\to\xi$. In view of (\ref{9}) and the Monotonicity formula applied to $\Sigma_k$ and $\xi_k$ the assumptions of Simon's graphical decomposition lemma (see \thref{Decomposition} in the appendix) are satisfied for $\rho\le\frac{\rho_0}{4}$ and infinitely many $k\in\varmathbb{N}$. Since $\Will(\Sigma_k)\le8\pi-\delta_0$, we can apply Lemma 1.4 in \cite{SL} to deduce that for $\theta\in\left(0,\frac{1}{2}\right)$ small enough, $\tau\in\left(\frac{\rho}{4},\frac{\rho}{2}\right)$ and infinitely many $k\in\varmathbb{N}$ only one of the discs $D_{\tau,l}^k$ appearing in the graphical decomposition lemma can intersect the ball $B_{\theta\frac{\rho}{4}}(\xi_k)$ (see \thref{Decomposition} for the notation). Moreover, by a slight perturbation from $\xi_k$ to $\xi$, we may assume that $\xi\in L_k$ for all $k\in\varmathbb{N}$. Now $L_k\to L$ in $\xi+G_2(\Rdrei)$, and therefore we may furthermore assume that the planes, on which the graph functions are defined, do not depend on $k\in\varmathbb{N}$. After all we get a graphical decomposition in the following way.

\begin{Le}\label{final}
For $\varepsilon\le\varepsilon_0$, $\rho\le\frac{\rho_0}{4}$ and infinitely many $k\in\varmathbb{N}$ there exist pairwise disjoint closed subsets $P_1^k,\ldots,P_{N_k}^k$ of $\Sigma_k$ such that 
$$\Sigma_k\cap\overline{B_{\theta\frac{\rho}{8}}(\xi)}=D_k\cap\overline{B_{\theta\frac{\rho}{8}}(\xi)}=\left(\graph u_k\cup\bigcup_n P_n^k\right)\cap\overline{B_{\theta\frac{\rho}{8}}(\xi)},$$
where $D_k$ is a topological disc and where the following holds:  
\begin{enumerate}
\item The sets $P_n^k$ are topological discs disjoint from $\graph u_k$.
\item $u_k\in C^\infty(\overline{\Omega_k},L^\perp)$, where $L\subset\Rdrei$ is a 2-dim. plane such that $\xi\in L$, and $\Omega_k=\left(B_{\lambda_k}(\xi)\cap L\right)\setminus\bigcup_m d_{k,m}$. Here $\lambda_k>\frac{\rho}{4}$ and the sets $d_{k,m}\subset L$ are pairwise disjoint closed discs.
\item The following inequalities hold:
\begin{eqnarray}
& & \hspace{-2cm}\sum_m\diam d_{k,m}+\sum_n\diam P_n^k\le c\left(\int_{\Sigma_k\cap B_{2\rho}(\xi)}|\rm{A}_k|^2\,d\Hzwei\right)^\frac{1}{4}\rho\le c\varepsilon^\frac{1}{2}\rho, \\
& & ||u_k||_{L^\infty(\Omega_k)}\le c\varepsilon^{\frac{1}{6}}\rho+\delta_k\quad\text{where }\delta_k\to0,\phantom{\left(\int_{B_{2\rho}}\right)^\frac{1}{4}} \\
& & ||\D u_k||_{L^\infty(\Omega_k)}\le c\varepsilon^{\frac{1}{6}}+\delta_k\quad\text{where }\delta_k\to0.\phantom{\left(\int_{B_{2\rho}}\right)^\frac{1}{4}}
\end{eqnarray}
\end{enumerate} 
\end{Le}

Now we leave the varifold context and define the functions
\begin{equation*}\label{12}
\schi=\chiup_{_{\Omega_k}}\in BV(\Rdrei),
\end{equation*}
where $\Omega_k\subset\Rdrei$ is the open, bounded set surrounded by $\Sigma_k$. We have
\begin{equation*}\label{13}
||\schi||_{L^1(\Rdrei)}=\frac{\sigma^3}{6\sqrt{\pi}}\quad\text{and}\quad|\D\schi|(U)=\mu_k(U)\le1\text{ for every open }U\subset\Rdrei.
\end{equation*}
Therefore the sequence $\schi$ is uniformly bounded in $BV(\Rdrei)$ and a compactness result for BV-functions (see \cite{EG}) yields that (after passing to a subsequence)
$$\schi\to\chiup\quad\text{in }L^1(\Rdrei)\text{ and pointwise a.e.}$$
for some function $\chiup\in BV(\Rdrei)$. Since the functions $\schi$ are characteristic functions we may assume without loss of generality that $\chiup$ is the characteristic function of a set $\Omega\subset\Rdrei$ with $\Ldrei(\Omega)=\frac{\sigma^3}{6\sqrt{\pi}}$. Because of the lower semicontinuity of the perimeter on open sets and the upper semicontinuity on compact sets under convergence of measures we get that
\begin{equation}\label{14}
|\D\chiup|\le\mu\quad\text{as measures}.
\end{equation}
In the end we would like to have that $\Sigma=\partial\Omega$ is smooth. Therefore it is necessary that $|\D\chiup|=\mu$ as measures, which actually holds.

\begin{Le}\label{chi=mu}
In the above setting we have for $\varepsilon\le\varepsilon_0$ that $|\D\chiup|=\mu$.
\end{Le}

\begin{Pro}
Let $\xi_0\in\Sigma\setminus\B$ and $\rho_0=\rho_0(\xi_0,\varepsilon)>0$ as in (\ref{9}). Let $\varepsilon\le\varepsilon_0$ such that \thref{final} holds and let $\rho\le\frac{\rho_0}{4}$. Let $\overline{u_k}\in C^{1,1}(B_{\lambda_k}(\xi_0)\cap L, L^\perp)$ be an extension of $u_k$ to the whole disc $B_{\lambda_k}(\xi_0)\cap L$ as in \thref{extension}, i.e. $\overline{u_k}=u_k$ in $\Omega_k$. From the $L^\infty$-bounds for the function $u_k$ and since $\lambda_k>\frac{\rho}{4}$ it follows that
\begin{eqnarray*}
\left\|\overline{u_k}\right\|_{L^\infty(B_\frac{\rho}{4}(\xi_0)\cap L)} & \le & c\varepsilon^\frac{1}{6}\rho+\delta_k\hspace{0,2cm}\le\hspace{0,2cm}c,\\
\left\|\D\overline{u_k}\right\|_{L^\infty(B_\frac{\rho}{4}(\xi_0)\cap L)} & \le & c\varepsilon^\frac{1}{6}+\delta_k\hspace{0,2cm}\le\hspace{0,2cm}c.
\end{eqnarray*}
Thus it follows that the sequence $\overline{u_k}$ is equicontinuous and uniformly bounded in $C^1(B_\frac{\rho}{4}(\xi_0)\cap L,L^\perp)$ and $W^{1,2}(B_\frac{\rho}{4}(\xi_0)\cap L,L^\perp)$. Therefore there exists a function $u\in C^{0,1}(B_\frac{\rho}{4}(\xi_0)\cap L,L^\perp)$ such that (after passing to a subsequence)
\begin{eqnarray*}
& & \hspace{-1cm}\overline{u_k}\to u\quad\text{in }C^0(B_\frac{\rho}{4}(\xi_0)\cap L,L^\perp),\phantom{\frac{1}{\rho}} \\
& & \hspace{-1cm}\overline{u_k}\rightharpoonup u\quad\text{weakly in }W^{1,2}(B_\frac{\rho}{4}(\xi_0)\cap L,L^\perp),\phantom{\frac{1}{\rho}} \\
& & \hspace{-1cm}\frac{1}{\rho}\left\|u\right\|_{L^\infty(B_\frac{\rho}{4}(\xi_0)\cap L)}+\left\|\D u\right\|_{L^\infty(B_\frac{\rho}{4}(\xi_0)\cap L)}\le c\varepsilon^\frac{1}{6}.
\end{eqnarray*}
Let $g\in C_c^1(B_{\theta\frac{\rho}{8}}(\xi_0),\Rdrei)$ with $|g|\le1$. It follows from the definition of $|\D\chiup|$ that
$$|\D\chiup|(B_{\theta\frac{\rho}{8}}(\xi_0))\ge\int\chiup\diver g=\lim_{k\to\infty}\int\schi\diver g.$$
\thref{final} yields
$$\int\schi\diver g=\int_{\graph u_k\cap B_{\theta\frac{\rho}{8}}(\xi_0)}\left< g,\nu_k\right>\,d\Hzwei+\sum_n\int_{P_n^k\cap B_{\theta\frac{\rho}{8}}(\xi_0)}\left< g,\nu_k\right>\,d\Hzwei,$$
where $\nu_k$ denotes the outer normal to $\partial\Omega_k=\Sigma_k$. Because of the Monotonicity formula and the diameter estimates for the sets $P_n^k$ we can estimate the second term on the right hand side by
$$\left|\sum_n\int_{P_n^k\cap B_{\theta\frac{\rho}{8}}(\xi_0)}\left< g,\nu_k\right>\,d\Hzwei\right|\le\sum_n\Hzwei(P_n^k)\le c\sum_n \left(\diam P_n^k\right)^2\le c\varepsilon\rho^2.$$
Because of the diameter estimates for the sets $d_{k,m}$ and the $L^\infty$-bounds for the functions $\overline{u_k}$ the first term on the right hand side can be estimated by
$$\left|\int_{\graph u_k\cap B_{\theta\frac{\rho}{8}}(\xi_0)}\left< g,\nu_k\right>\,d\Hzwei\right|\ge\left|\int_{\graph\overline{u_k}\cap B_{\theta\frac{\rho}{8}}(\xi_0)}\left< g,\overline{\nu}_k\right>\,d\Hzwei\right|-c\varepsilon\rho^2,$$
where $\overline{\nu}_k$ denotes the outer normal to $\graph\overline{u_k}$. Using the convergence stated above together with the estimates for the limit function $u$ we get that 
$$\liminf_{k\to\infty}\left|\int_{\graph\overline{u_k}\cap B_{\theta\frac{\rho}{8}}(\xi_0)}\left< g,\overline{\nu}_k\right>\,d\Hzwei\right|\ge\left(\frac{\theta}{8}-c\varepsilon^\frac{1}{6}\right)^2\pi\rho^2$$
and therefore
\begin{equation*}
|\D\chiup|(B_{\theta\frac{\rho}{8}}(\xi_0))\ge\left(\frac{\theta}{8}-c\varepsilon^\frac{1}{6}\right)^2\pi\rho^2-c\varepsilon\rho^2.
\end{equation*}
In the same way, using that $\mu_k\to\mu$ in $C_c^0(\Rdrei)'$, we get
\begin{equation*}
\mu(B_{\theta\frac{\rho}{8}}(\xi_0))\le\sqrt{1+c\varepsilon^\frac{1}{3}}\left(\frac{\theta}{8}+c\varepsilon^\frac{1}{6}\right)^2\pi\rho^2+c\varepsilon\rho^2.
\end{equation*}
Since the derivation $\D_\mu|\D\chiup|(\xi_0)$ exists for $\mu$-a.e. $\xi_0\in\Sigma$ (see \cite{EG}) we get for $\mu$-a.e. $\xi_0\in\Sigma\setminus\B$ that
$$\D_\mu|\D\chiup|(\xi_0)\ge\frac{\left(1-c\varepsilon^\frac{1}{6}\right)^2\pi-c\varepsilon}{\sqrt{1+c\varepsilon^\frac{1}{3}}\left(1+c\varepsilon^\frac{1}{6}\right)^2\pi+c\varepsilon}.$$
Letting $\varepsilon=\frac{1}{n}\to0$ we get in view of (\ref{14})
\begin{equation*}
\D_\mu|\D\chiup|(\xi_0)=1\quad\text{for }\mu\text{-a.e. }\xi_0\in\Sigma\setminus\bigcup_n\mathcal{B}_\frac{1}{n}.
\end{equation*}
Since each set $\mathcal{B}_\frac{1}{n}$ contains only finitely many points and since the Monotonicity formula yields $\mu(\{\xi\})=0$ for every $\xi\in\Rdrei$, we get that
$$\D_\mu|\D\chiup|(\xi_0)=1\quad\text{for }\mu\text{-a.e. }\xi_0\in\Rdrei.$$
Using the theorem of Radon-Nikodym the lemma follows from (\ref{14}).
\end{Pro}

\begin{Rem}
\rm{Notice that the only thing we needed up to now was the bound on the Willmore energy $\Will(\Sigma_k)\le8\pi-\delta_0$. We are now able to prove that 
\begin{equation}\label{300}
\lim_{\sigma\searrow0}\beta(\sigma)=8\pi.
\end{equation}
We already know that $\beta$ is decreasing and bounded by $8\pi$. Therefore the limit exists. Let $\sigma_l\to0$ and assume (\ref{300}) is false. After passing to a subsequence there exists a $\delta_0>0$ such that $\beta(\sigma_l)\le8\pi-\delta_0$ for all $l\in\varmathbb{N}$. Let $\Sigma_l\in\mathcal{M}_{\sigma_l}$ such that $\Will(\Sigma_l)\le\beta(\sigma_l)+\frac{\delta_0}{2}\le8\pi-\frac{\delta_0}{2}$ and let $\Omega_l\subset\Rdrei$ be the open set surrounded by $\Sigma_l$. Again after scaling and translation we may assume that $\Hzwei(\Sigma_l)=1$ and $0\in\Sigma_l$, and that the radon measures $\mu_l=\Hzwei\llcorner\Sigma_l$ converge to a radon measure $\mu$ with $\mu(\Rdrei)=1$. On the other hand we have that the $BV$-functions $\chiup_l=\chiup_{\Omega_l}$ are uniformly bounded and therefore converge (after passing to a subsequence) in $L^1$ to a $BV$-function $\chiup$. Since $I(\Sigma_l)\to0$ and $\Hzwei(\Sigma_l)=1$ it follows that $\chiup=0$. Finally, since $\Will(\Sigma_l)\le8\pi-\frac{\delta_0}{2}$, we can do exactly the same as before to get $\mu=|\D\chiup|$, which contradicts $\mu(\Rdrei)=1$. Therefore (\ref{300}) holds.}
\end{Rem}

We continue with the proof of \thref{MainThm}. The main idea to prove regularity is to derive a power decay for the $L^2$-norm of the second fundamental form via constructing comparison surfaces by a cut-and-paste procedure as done in \cite{SL}. But this method cannot be directly applied in our case, since the isoperimetric ratio might change by this procedure. In order to correct the isoperimetric ratio of the generated surfaces, we will apply an appropriate variation. But what is an appropriate variation in our case? Which is the quantity we have to look at? To answer this question let $\Phi:(-\varepsilon,\varepsilon)\times\Rdrei\to\Rdrei$ be a $C^2$-variation with compact support and define $\Omega_{k,t}=\Phi_t(\Omega_k)$, $\Sigma_{k,t}=\partial\Omega_{k,t}=\Phi_t(\Sigma_k)$ and $X(x)=\partial_t\Phi_t(x)_{|_{t=0}}$. It follows that
\begin{equation}\label{1000}
\frac{d}{dt}I(\Sigma_{k,t})_{|_{t=0}}=\frac{I(\Sigma_k)}{3\Hzwei(\Sigma_k)}\left(\frac{3}{2}\int\big<X,\vec{H}_k\big>\,d\mu_k+\frac{\Hzwei(\Sigma_k)}{\Ldrei(\Omega_k)}\int\chiup_{\Omega_k}\diver X\right).
\end{equation}
Because of (\ref{4''}) and since $\chiup_{\Omega_k}\to\chiup_\Omega$, we get in view of \thref{chi=mu} that
\begin{equation}\label{3003}
\lim_{k\to\infty}\frac{d}{dt}I(\Sigma_{k,t})_{|_{t=0}}=\frac{\sigma}{3}\int\left<X,\frac{3}{2}\vec{H}-\frac{6\sqrt{\pi}}{\sigma^3}\nu\right>\,d\mu,
\end{equation}
where $\nu$ is given by the equation 
$$\int\chiup\diver g=-\int\big<g,\nu\big>\,d|\D\chiup|=-\int\big<g,\nu\big>\,d\mu$$
for $g\in C_c^1(\Rdrei,\Rdrei)$. This follows from the Riesz representation theorem applied to BV-functions and \thref{chi=mu}.\\[0,3cm]
Now if there would exist a vector field $X\in C^\infty_c(\Rdrei,\Rdrei)$ such that the right hand side of (\ref{3003}) is not equal to $0$, we would have that the first variation of the isoperimetric ratio of $\Sigma_k$ not equals $0$ for $k$ sufficiently large, and in conclusion we would have a chance to correct the isoperimetric ratio of the generated surfaces. The next lemma is concerned with the existence of such a vector field and relies on the fact that each surface $\Sigma\in\M$ is not a round sphere.

\begin{Le}\label{vectorfield}
There exists a $R>0$ such that for every $\xi\in\Sigma$ there exists a point $\eta\in\Sigma\setminus B_R(\xi)$ such that for all $\beta>0$ there exists a vectorfield $X\in C_c^\infty(B_\beta(\eta),\Rdrei)$ such that
$$\int\left<X,\frac{3}{2}\vec{H}-\frac{6\sqrt{\pi}}{\sigma^3}\nu\right>\,d\mu\neq0.$$
\end{Le}

\begin{Pro}
Assume the statement is false. Then there exists a sequence $R_k\searrow0$ and $\xi_k\in\Sigma$ such that for all $\eta\in\Sigma\setminus B_{R_k}(\xi_k)$ there exists a $\beta_\eta>0$ such that 
$$\int\left<X,\frac{3}{2}\vec{H}-\frac{6\sqrt{\pi}}{\sigma^3}\nu\right>\,d\mu=0$$
for all $X\in C_c^\infty(B_{\beta_\eta}(\eta),\Rdrei)$. Since $\Sigma$ is compact, it follows after passing to a subsequence that $\xi_k\to\xi\in\Sigma$, and since $\mu(\{\xi\})=0$, which follows from \thref{Monotonie}, we get after all that
\begin{equation}
\vec{H}(x)=\frac{4\sqrt{\pi}}{\sigma^3}\nu(x)\quad\text{for }\mu\text{-a.e. }x\in\Sigma.
\end{equation}
Now the idea of the proof is the following: We just have to show that $\Sigma$ is smooth, because then $\Sigma$ would be a smooth surface with constant mean curvature and Willmore energy smaller than $8\pi$. By a theorem of Alexandroff $\Sigma$ would be a round sphere which contradicts our choice of $\sigma\in(0,1)$. To show that $\Sigma$ is smooth we just have to show that $\theta^2(\mu,x)=1$ for every $x\in\Sigma$, because then Allard's regularity theorem would yield (remember that $\vec{H}\in L^\infty(\mu)$ now) that $\Sigma$ can be written as a $C^{1,\alpha}$-graph around $x$ that solves the constant mean curvature equation and is therefore smooth.\\[0,3cm]
Let $x_0\in\Sigma$. To prove that $\theta^2(\mu,x_0)=1$ notice that, since $\chiup\in BV(\Rdrei)$, $\mu$ generates an integer multiplicity, rectifiable 2-current $\cur_\mu\in\mathcal{D}^2(\Rdrei)'$ with $\partial \cur_\mu=0$. Denote by $\mu_{x_0,\lambda}$ the blow-ups of $\mu$ around $x_0$. Now also the blow-ups generate integer multiplicity, rectifiable 2-currents $\cur_{x_0,\lambda}$ with $\partial\cur_{x_0,\lambda}=0$. Moreover the mass of $\cur_{x_0,\lambda}$ of a set $W\subset\subset\Rdrei$ such that $W\subset B_R(0)$ is estimated in view of the Monotonicity formula by
\begin{eqnarray*}
\underline{M}_W(\cur_{x_0,\lambda})\le\mu_{x_0,\lambda}(B_R(0))=\lambda^{-2}\mu(B_{\lambda R}(x_0))\le cR^2.
\end{eqnarray*}   
By a compactness theorem for integer multiplicity, rectifiable 2-currents (see \cite{SLGMT}) there exists an integer multiplicity, rectifiable 2-current $\cur_{x_0}\in\mathcal{D}^2(\Rdrei)'$ such that $\partial\cur_{x_0}=0$ and (after passing to a subsequence)
$$\cur_{x_0,\lambda}\to\cur_{x_0}\quad\text{for }\lambda\to0\text{ weakly as currents}.$$
Let $\mu_{x_0}$ be the underlying varifold.\\[0,3cm]
On the other hand there exists a stationary, integer multiplicity, rectifiable 2-cone $\mu_\infty$ such that (after passing to a subsequence)
$$\mu_{x_0,\lambda}\to\mu_\infty\quad\text{for }\lambda\to0\text{ weakly as varifolds}.$$
Now we get the following:
\begin{itemize}
\item[1.)] $\mu_{x_0}\le\mu_\infty$: This follows from the lower semicontinuity of the mass with respect to weak convergence of currents and the upper semicontinuity on compact sets with respect to weak convergence of measures.
\item[2.)] $\theta^2(\mu_\infty,\cdot)\le2-\frac{\delta_0}{4\pi}$ everywhere: Since $\mu_\infty$ is a stationary 2-cone, the Monotonicity formula yields for all $z\in\Rdrei$ and all $B_\tau(0)$ such that $\mu_\infty(\partial B_\tau(0))=0$
$$\theta^2(\mu_\infty,z)\le\theta^2(\mu_\infty,0)=\frac{\mu_\infty(B_\tau(0))}{\pi\tau^2}=\liminf_{\lambda\to0}\frac{\mu_{x_0,\lambda}(B_\tau(0))}{\pi\tau^2}.$$
Now since $\frac{\mu_{x_0,\lambda}(B_\tau(0))}{\pi\tau^2}=\frac{\mu(B_{\lambda\tau}(x_0))}{\pi(\lambda\tau)^2}$, it follows that $\theta^2(\mu_\infty,z)\le\underline{\theta}\,^2(\mu,x_0)$, and the claim follows from (\ref{4}).
\item[3.)] $\theta^2(\mu_\infty,\cdot)=1$ $\mu_\infty$-a.e.: This follows from 2.) since $\mu_\infty$ is integral.
\item[4.)] $\mu_{x_0}=\mu_\infty$: Choose a point $x\in\Rdrei$ such that $\theta^2(\mu_\infty,x)=1$. By Allard's regularity theorem there exists a neighborhood $U(x)$ of $x$ in which $\mu_\infty$ can be written as a $C^{1,\alpha}$-graph, which is actually smooth since $\mu_\infty$ is stationary. Moreover we get that the convergence $\mu_{x_0,\lambda}\llcorner U(x)\to\mu_\infty\llcorner U(x)$ is in $C^{1,\alpha}$. Thus $\mu_\infty\llcorner U(x)\leftarrow\mu_{x_0,\lambda}\llcorner U(x)\rightarrow\mu_{x_0}\llcorner U(x)$, hence for $\mathcal{U}=\bigcup_{\theta^2(\mu_\infty,x)=1}U(x)$ we have that $\mu_\infty\llcorner\mathcal{U}=\mu_{x_0}\llcorner\mathcal{U}$. Since we already know that $\theta^2(\mu_\infty,x)=1$ for $\mu_\infty$-a.e. $x\in\Rdrei$ we get 4.).
\end{itemize}
From 4.) it follows that $\cur_{x_0}$ is a stationary, integer multiplicity, rectifiable 2-current with $\partial\cur_{x_0}=0$. Since moreover $\mu_{x_0}=\mu_\infty$ is a stationary, rectifiable 2-cone we get for all $\tau>0$ and all $z\in\Rdrei$ that
$$\frac{\mu_{x_0}(B_\tau(z))}{\pi\tau^2}\le\theta^2(\mu_{x_0},0)\le2-\frac{\delta_0}{4\pi}.$$
Letting $\tau\to\infty$ we get that $\theta^2(\mu_{x_0},\infty)\le2-\frac{\delta_0}{4\pi}$. Using Theorem 2.1 in \cite{KLS} it follows that $\cur_{x_0}$ is a unit density plane or
$$\mu_{x_0}=\mu_\infty=\Hzwei\llcorner P\quad\text{for some }P\in G_2(\Rdrei).$$
Therefore we get for all balls $B_\tau(0)$ such that $\mu_\infty(\partial B_\tau(0))=0$
$$\theta^2(\mu,x_0)=\lim_{\lambda\to0}\frac{\mu(B_{\lambda\tau}(x_0))}{\pi(\lambda\tau)^2}=\lim_{\lambda\to0}\frac{\mu_{x_0,\lambda}(B_{\tau}(0))}{\pi\tau^2}=\frac{\mu_\infty(B_\tau(0))}{\pi\tau^2}=\frac{\Hzwei\llcorner P(B_\tau(0))}{\pi\tau^2}=1$$
and the lemma is proved.
\end{Pro}

In the next step we prove a power decay for the $L^2$-norm of the second fundamental form on small balls around the good points $\xi\in\Sigma\setminus\B$. This will help us to show that $\Sigma$ is actually $C^{1,\alpha}\cap W^{2,2}$ away from the bad points.

\begin{Le}\label{2ff-absch}
Let $\xi_0\in\Sigma\setminus\B$. There exists a $\rho_0=\rho_0(\xi_0,\varepsilon)>0$ such that for all $\xi\in\Sigma\cap B_{\frac{\rho_0}{2}}(\xi_0)$ and all $\rho\le\frac{\rho_0}{4}$ we have
$$\liminf_{k\to\infty}\int_{\Sigma_k\cap B_{\theta\frac{\rho}{8}}(\xi)}|\rm{A}_k|^2\,d\Hzwei\le c\rho^\alpha,$$
where $\alpha\in(0,1)$ and $c<\infty$ are universal constants.
\end{Le}

\begin{Pro}
Choose according to \thref{vectorfield} a $R>0$ such that for every $\xi\in\Sigma$ there exists a point $\eta\in\Sigma\setminus B_R(\xi)$ such that for all $\beta>0$ there exists a vectorfield $X\in C_c^\infty(B_\beta(\eta),\Rdrei)$ such that
\begin{equation}\label{2000}
\int\left<X,\frac{3}{2}\vec{H}-\frac{6\sqrt{\pi}}{\sigma^3}\nu\right>\,d\mu\neq0.
\end{equation}
Let $\xi_0\in\Sigma\setminus\B$, $\rho_0>0$ as in (\ref{9}). We may assume without loss of generality that 
\begin{equation}
\rho_0<\frac{R}{2}.
\end{equation}
Let $\xi\in\Sigma\cap B_{\frac{\rho_0}{2}}(\xi_0)$ and $\rho\le\frac{\rho_0}{4}$. Notice that \thref{final} holds. For $\tau\in\left(\theta\frac{\rho}{16},\frac{3}{4}\theta\frac{\rho}{8}\right)$ define the set
$$C_\tau(\xi)=\left\{x+y\,\big|\,x\in B_\tau(\xi)\cap L, y\in L^\perp\right\}.$$
From the $L^\infty$-estimates for the functions $u_k$ and the diameter estimates for the sets $P_n^k$ it follows for $\varepsilon\le\varepsilon_0$ and $\delta_k\le\frac{1}{8}\theta\frac{\rho}{8}$ that $D_k\cap C_\tau(\xi)=D_k\cap C_\tau(\xi)\cap\overline{B_{\theta\frac{\rho}{8}}(\xi)}$. Therefore
$$\Sigma_k\setminus\Big(D_k\cap C_\tau(\xi)\Big)=\Sigma_k\setminus\left(C_\tau(\xi)\cap\overline{B_{\theta\frac{\rho}{8}}(\xi)}\right)\quad\text{for }\varepsilon\le\varepsilon_0\text{ and }\delta_k\le\frac{1}{8}\theta\frac{\rho}{8}.$$
Define the sets
$$S_k(\xi)=\left\{\tau\in\left(\theta\frac{\rho}{16},\frac{3}{4}\theta\frac{\rho}{8}\right)\,\Bigg|\,\partial C_\tau(\xi)\cap\bigcup_m d_{k,m}=\emptyset\right\},$$
$$T_k(\xi)=\left\{\tau\in S_k(\xi)\,\Bigg|\,\int_{D_k\cap\partial C_\tau(\xi)}|\rm{A}_k|^2\,d\Hzwei\le\frac{128}{\theta\rho}\int_{D_k\cap C_{\frac{3}{4}\theta\frac{\rho}{8}}(\xi)\setminus C_{\theta\frac{\rho}{16}}(\xi)}|\rm{A}_k|^2\,d\Hzwei\right\}.$$
Using the diameter estimates for the discs $d_{k,m}$ we get that $\Leins(S_k(\xi))\ge\theta\frac{\rho}{64}$ for $\varepsilon\le\varepsilon_0$, and then from a Fubini-type argument that $\Leins(T_k(\xi))\ge\theta\frac{\rho}{128}$. From the selection principle in \cite{SL}, Lemma B.1, it follows that there exists a $\tau\in\left(\theta\frac{\rho}{16},\frac{3}{4}\theta\frac{\rho}{8}\right)$ such that $\tau\in T_k(\xi)$ for infinitely many $k\in\varmathbb{N}$. \\[0,3cm]
Apply \thref{extension} to get a function $w_k\in C^\infty(\overline{B_\tau(\xi)}\cap L, L^\perp)$ for infinitely many $k\in\varmathbb{N}$ such that
\begin{eqnarray*}
& (i) & w_k=u_k\quad\text{and}\quad\frac{\partial w_k}{\partial\nu}=\frac{\partial u_k}{\partial\nu}\quad\text{on }\partial B_\tau(\xi)\cap L,\phantom{\int_{d_{k,m}^\sim}}\\
& (ii) & \frac{1}{\tau}||w_k||_{L^\infty(B_\tau(\xi)\cap L)}\le c\varepsilon^\frac{1}{6}+\frac{\delta_k}{\tau}\quad\text{where }\delta_k\to0,\phantom{\int_{d_{k,m}^\sim}}\\
& (iii) & ||\D w_k||_{L^\infty(\partial B_\tau(\xi)\cap L)}\le c\varepsilon^{\frac{1}{6}}+\delta_k\quad\text{where }\delta_k\to0,\phantom{\int_{d_{k,m}^\sim}}\\
& (iv) & \int_{B_\tau(\xi)\cap L}|\D^2 w_k|^2\le c\tau\int_{\graph {u_k}_{|_{\partial B_\tau(\xi)\cap L}}}|\rm{A}_k|^2\,d\Heins.
\end{eqnarray*}
Since $\graph w_k\subset\overline{B_{\theta\frac{\rho}{8}}(\xi)}$ for $\varepsilon\le\varepsilon_0$ and $\delta_k\le\frac{1}{8}\theta\frac{\rho}{8}$ we get from the above that
$$\graph w_k\cap\Big(\Sigma_k\setminus\Big(D_k\cap C_\tau(\xi)\Big)\Big)\subset C_\tau(\xi)\cap\overline{B_{\theta\frac{\rho}{8}}(\xi)}\cap\left(\Sigma_k\setminus\left(C_\tau(\xi)\cap\overline{B_{\theta\frac{\rho}{8}}(\xi)}\right)\right)=\emptyset.$$
Now define the surfaces
$$\tilde\Sigma_k=\Sigma_k\setminus\Big(D_k\cap C_\tau(\xi)\Big)\cup\graph w_k.$$
From the above it follows for $\varepsilon\le\varepsilon_0$ and $\delta_k\le\frac{1}{8}\theta\frac{\rho}{8}$ that $\tilde\Sigma_k$ is a compact, embedded and connected $C^{1,1}$-surface with $\genus\tilde\Sigma_k=0$. In addition $\tilde\Sigma_k$ surrounds an open set $\tilde\Omega_k$ and $\tilde\Sigma_k\cap B_\frac{R}{2}(\eta)=\Sigma_k\cap B_\frac{R}{2}(\eta)$.\\[0,3cm]
Next we compare the isoperimetric coefficients of $\Sigma_k$ and $\tilde\Sigma_k$. Using the $L^\infty$-bounds for $w_k$ and the Monotonicity formula we get from the definition of $\tilde\Sigma_k$ that
\begin{equation}\label{600}
\left|\Hzwei(\tilde\Sigma_k)-\Hzwei\left(\Sigma_k\right)\right|\le\Hzwei(\Sigma_k\cap B_{\theta\frac{\rho}{8}}(\xi))+\Hzwei(\graph w_k)\le c\rho^2.
\end{equation}
Since $|\Ldrei(\tilde\Omega_k)-\Ldrei(\Omega_k)|\le\Ldrei(\Omega\Delta\tilde\Omega)$ and since by construction $\Omega\Delta\tilde\Omega\subset B_{\theta\frac{\rho}{8}}(\xi)$, it follows that
\begin{equation}\label{602}
\left|\Ldrei(\tilde\Omega_k)-\Ldrei(\Omega_k)\right|\le c\rho^3.
\end{equation}
Since $\Hzwei(\Sigma_k)=1$ and $\Ldrei(\Omega_k)=\frac{\sigma^3}{6\sqrt{\pi}}$, we get by choosing $\rho_0$ smaller (smaller in an universal way) that
\begin{equation}\label{601}
\frac{1}{2}\le\Hzwei(\tilde\Sigma_k)\le2,
\end{equation}
\begin{equation}\label{603}
\frac{\sigma^3}{12\sqrt{\pi}}\le\Ldrei(\tilde\Omega_k)\le\frac{\sigma^3}{3\sqrt{\pi}}.
\end{equation}
Moreover we finally get 
\begin{equation}\label{18}
\left|I(\tilde\Sigma_k)-\sigma\right|=\left|I(\tilde\Sigma_k)-I(\Sigma_k)\right|\le c\rho,
\end{equation}
and we may assume without loss of generality that
\begin{equation}\label{18'}
\frac{\sigma}{2}\le I(\tilde\Sigma_k)\le2\sigma.
\end{equation}
As mentioned before $\tilde\Sigma_k$ might not have the right isoperimetric ratio and may therefore not be a comparison surface.\\[0,3cm] 
According to (\ref{2000}) let $\eta\in\Sigma\setminus B_R(\xi)$ and $X\in C_c^\infty(B_\frac{R}{2}(\eta),\Rdrei)$ such that
\begin{equation}\label{2001}
\int\left<X,\frac{3}{2}\vec{H}-\frac{6\sqrt{\pi}}{\sigma^3}\nu\right>\,d\mu\ge c_0>0.
\end{equation}
Notice that the constant $c_0$ does not depend on $\varepsilon$, $\xi$, $\rho$ or $k\in\varmathbb{N}$.\\[0,3cm]
Let $\Phi\in C^\infty(\Reins\times\Rdrei,\Rdrei)$ be the flow of the vectorfield $X$, namely
\begin{eqnarray*}
& & \Phi_t(\cdot)=\Phi(t,\cdot)\in C^\infty(\Rdrei,\Rdrei)\quad\text{is a diffeomorphism for all }t\in\Reins, \\
& & \Phi(0,z)=z\quad\text{for all }z\in\Rdrei, \\
& & \partial_t\Phi(t,z)=X(\Phi(t,z))\quad\text{for all }(t,z)\in\Reins\times\Rdrei.
\end{eqnarray*}
Since $\spt X\subset B_\frac{R}{2}(\eta)$ there exists a $T_0=T_0(X)>0$ such that for all $t\in(-T_0,T_0)$
$$\Phi_t=Id\quad\text{in }\Rdrei\setminus B_\frac{R}{2}(\eta).$$
Define the sets
\begin{equation}
\tilde\Omega_k^t=\Phi_t(\tilde\Omega_k)\quad\text{and}\quad\tilde\Sigma_k^t=\partial\tilde\Omega_k^t=\Phi_t(\tilde\Sigma_k).
\end{equation}
Choosing $T_0$ smaller if necessary (depending on $X$) it follows for $t\in(-T_0,T_0)$ that
$$\Hzwei(\tilde\Sigma_k^t)=\int_{\tilde\Sigma_k}J_{\tilde\Sigma_k}\Phi_t\,d\Hzwei\quad\text{and}\quad\Ldrei(\tilde\Omega_k^t)=\int_{\tilde\Omega_k}\det\D\Phi_t.$$
By choosing $T_0$ smaller if necessary (depending on $X$) and estimating very roughly we get that there exists a constant $0<c=c(X)<\infty$ such that for all $t\in(-T_0,T_0)$
\begin{eqnarray*}
& (i) & \frac{1}{c}\Hzwei(\tilde\Sigma_k)\le\sup_{t\in(-T_0,T_0)}\Hzwei(\tilde\Sigma_k^t)\le c\Hzwei(\tilde\Sigma_k),\phantom{\int_{\tilde\Sigma_k^t}|\rm{A}_k^t|^2} \\
& (ii) & \frac{1}{c}\Ldrei(\tilde\Omega_k)\le\sup_{t\in(-T_0,T_0)}\Ldrei(\tilde\Omega_k^t)\le c\Ldrei(\tilde\Omega_k),\phantom{\int_{\tilde\Sigma_k^t}|\rm{A}_k^t|^2} \\
& (iii) & \sup_{t\in(-T_0,T_0)}\left|\frac{d}{dt}\Hzwei(\tilde\Sigma_k^t)\right|+\sup_{t\in(-T_0,T_0)}\left|\frac{d}{dt}\Ldrei(\tilde\Omega_k^t)\right|\le c,\phantom{\int_{\tilde\Sigma_k^t}|\rm{A}_k^t|^2} \\
& (iv) & \sup_{t\in(-T_0,T_0)}\left|\frac{d^2}{dt^2}\Hzwei(\tilde\Sigma_k^t)\right|+\sup_{t\in(-T_0,T_0)}\left|\frac{d^2}{dt^2}\Ldrei(\tilde\Omega_k^t)\right|\le c,\phantom{\int_{\tilde\Sigma_k^t}|\rm{A}_k^t|^2} \\
& (v) & \sup_{t\in(-T_0,T_0)}\left|\frac{d}{dt}\int_{\tilde\Sigma_k^t}|\rm{A}_k^t|^2\,d\Hzwei\right|\le c.
\end{eqnarray*}
The last inequality can be proved by writing $\tilde\Sigma_k$ locally as a graph with small Lipschitz norm and using a partition of unity. \\[0,3cm]
Now first of all it follows for the first variation of the isoperimetric coefficient of $\tilde\Sigma_k$, using that $\spt X\subset B_\frac{R}{2}(\eta)$ and $\tilde\Sigma_k\cap B_\frac{R}{2}(\eta)=\Sigma_k\cap B_\frac{R}{2}(\eta)$,
$$\frac{d}{dt}I(\tilde\Sigma_k^t)_{|_{t=0}}=\frac{I(\tilde\Sigma_k)}{3\Hzwei(\tilde\Sigma_k)}\int_{\Sigma_k}\Bigg<X,\frac{3}{2}\vec{H}_k-\frac{\Hzwei(\tilde\Sigma_k)}{\Ldrei(\tilde\Omega_k)}\nu_k\Bigg>\,d\Hzwei.$$
Now it follows from (\ref{600})-(\ref{603}) that
$$\left|\int_{\Sigma_k}\Bigg<X,\left(\frac{\Hzwei(\tilde\Sigma_k)}{\Ldrei(\tilde\Omega_k)}-\frac{6\sqrt{\pi}}{\sigma^3}\right)\nu_k\Bigg>\,d\Hzwei\right|\le c\left|\frac{\Hzwei(\tilde\Sigma_k)}{\Ldrei(\tilde\Omega_k)}-\frac{\Hzwei(\Sigma_k)}{\Ldrei(\Omega_k)}\right|\Hzwei(\Sigma_k)\le c\rho,$$
where $c=c(X)$, and therefore (\ref{601}) and (\ref{18'}) yield
$$\frac{d}{dt}I(\tilde\Sigma_k^t)_{|_{t=0}}\ge\frac{I(\tilde\Sigma_k)}{3\Hzwei(\tilde\Sigma_k)}\int_{\Sigma_k}\Bigg<X,\frac{3}{2}\vec{H}_k-\frac{6\sqrt{\pi}}{\sigma^3}\nu_k\Bigg>\,d\Hzwei-c\rho.$$
Since $\int_{\Sigma_k}\Big<X,\frac{3}{2}\vec{H}_k-\frac{6\sqrt{\pi}}{\sigma^3}\nu_k\Big>\,d\Hzwei\to\int\Big<X,\frac{3}{2}\vec{H}-\frac{6\sqrt{\pi}}{\sigma^3}\nu\Big>\,d\mu\ge c_0$, it follows from (\ref{601}) and (\ref{18'}) that there exists a constant $0<c_0<\infty$ independent of $\varepsilon$, $\xi$, $\rho$ and $k\in\varmathbb{N}$, such that for $k$ sufficiently large
\begin{equation}\label{3000}
\frac{d}{dt}I(\tilde\Sigma_k^t)_{|_{t=0}}\ge c_0-c\rho.
\end{equation}
Moreover using the estimates (iii) and (iv) above it follows that
\begin{equation}
\sup_{t\in(-T_0,T_0)}\left|\frac{d^2}{dt^2}I(\tilde\Sigma_k^t)\right|\le c,
\end{equation}
where $c=c(X)<\infty$ is a universal constant.\\[0,3cm]
Using Taylor's formula we get in view of (\ref{18}) that for each $k\in\varmathbb{N}$ there exists a $t_k$ with $|t_k|\le c\rho$ such that
$$I(\tilde\Sigma_k^{t_k})=\sigma.$$
Therefore we get by construction that $\tilde\Sigma_k^{t_k}\in\M$ is a comparison surface to $\Sigma_k$. Moreover it follows from (v) above that
$$\left|\int_{\tilde\Sigma_k^{t_k}}|\rm{A}_k^{t_k}|^2\,d\Hzwei-\int_{\tilde\Sigma_k}|\rm{A}_k|^2\,d\Hzwei\right|\le\left|t_k\right|\sup_{t\in[-t_k,t_k]}\left|\frac{d}{dt}\int_{\tilde\Sigma_k^t}|\rm{A}_k^t|^2\,d\Hzwei\right|\le c\rho.$$
Since $\Sigma_k$ is a minimizing sequence for the Willmore functional in $\M$ and by the Gauss-Bonnet theorem therefore a minimizing sequence for the functional $\int_{\Sigma}|\rm{A}|^2$, we get
\begin{eqnarray*}
\int_{\Sigma_k}|\rm{A}_k|^2\,d\Hzwei\le\int_{\tilde\Sigma_k}|\rm{A}_k|^2\,d\Hzwei+c\rho+\varepsilon_k\quad\text{with }\varepsilon_k\to0.
\end{eqnarray*}
Now by definition of $\tilde\Sigma_k$ it follows that
$$\int_{D_k\cap C_\tau(\xi)}|\rm{A}_k|^2\,d\Hzwei\le\int_{\graph w_k}|\rm{A}_k|^2\,d\Hzwei+c\rho+\varepsilon_k.$$
By definition of $w_k$ and the choice of $\tau$ we get 
$$\int_{\graph w_k}|\rm{A}_k|^2\,d\Hzwei\le c\int_{D_k\cap C_{\frac{3}{4}\theta\frac{\rho}{8}}(\xi)\setminus C_{\theta\frac{\rho}{16}}(\xi)}|\rm{A}_k|^2\,\Hzwei.$$
Since $B_{\theta\frac{\rho}{16}}(\xi)\subset C_\tau(\xi)$, we get that (remember that $D_k\cap B_{\theta\frac{\rho}{8}}(\xi)=\Sigma_k\cap B_{\theta\frac{\rho}{8}}(\xi)$)
$$\int_{\Sigma_k\cap B_{\theta\frac{\rho}{16}}(\xi)}|\rm{A}_k|^2\,d\Hzwei\le c\int_{\Sigma_k\cap B_{\theta\frac{\rho}{8}}(\xi)\setminus B_{\theta\frac{\rho}{16}}(\xi)}|\rm{A}_k|^2\,d\Hzwei+c\rho+\varepsilon_k.$$
Now by adding $c$ times the left hand side of this inequality to both sides ("hole filling") we deduce the following:\\[0,3cm]
For $\rho\le\frac{\rho_0}{4}$ and infinitely many $k\in\varmathbb{N}$ it follows that
$$\int_{\Sigma_k\cap B_{\theta\frac{\rho}{16}}(\xi)}|\rm{A}_k|^2\,d\Hzwei\le\gamma\int_{\Sigma_k\cap B_{\theta\frac{\rho}{8}}(\xi)}|\rm{A}_k|^2\,d\Hzwei+c\rho+\varepsilon_k.$$
where $\gamma=\frac{c}{c+1}\in(0,1)$ is a fixed universal constant. If we let 
$$g(\rho)=\liminf_{k\to\infty}\int_{\Sigma_k\cap B_{\theta\frac{\rho}{16}}(\xi)}|\rm{A}_k|^2\,d\Hzwei$$
we get that
$$g(\rho)\le\gamma g(2\rho)+c\rho\quad\text{for all }\rho\le\frac{\rho_0}{4}.$$
Now in view of \thref{decay} it follows that
$$g(\rho)\le c\rho^\alpha\quad\text{for all }\rho\le\frac{\rho_0}{2}$$
and the lemma is proved.
\end{Pro}

In the next step we want to do the same as in the proof of \thref{chi=mu} where we constructed a sequence of functions which converged strongly in $C^0$ and weakly in $W^{1,2}$. But now with the estimate of \thref{2ff-absch} we will get better control on the sequence.\\[0,3cm]
So let $\xi\in\Sigma\cap B_{\frac{\rho_0}{2}}(\xi_0)$. Define the quantity $\alpha_k(\rho)$ by
$$\alpha_k(\rho)=\int_{\Sigma_k\cap B_{2\rho}(\xi)}|\rm{A}_k|^2\,d\Hzwei$$
and notice that by the choice of $\rho_0$ and \thref{2ff-absch} we have that
\begin{equation}\label{20}
\alpha_k(\rho)\le c\varepsilon^2\quad\text{and}\quad\liminf_{k\to\infty}\alpha_k(\rho)\le c\rho^\alpha\quad\text{for all }\rho\le\theta\frac{\rho_0}{64}.
\end{equation}
Furthermore we get from \thref{final} and the Monotonicity formula that
\begin{equation}\label{21}
\sum_m\diam d_{k,m}\le c\alpha_k(\rho)^\frac{1}{4}\rho\le c\varepsilon^\frac{1}{2}\rho\quad\text{and}\quad\sum_m\Lzwei\left(d_{k,m}\right)\le c\alpha_k(\rho)^\frac{1}{2}\rho^2,
\end{equation}
\begin{equation}\label{21'}
\sum_n\diam P_n^k\le c\alpha_k(\rho)^\frac{1}{4}\rho\quad\text{and}\quad\sum_n\Hzwei\left(P_n^k\right)\le c\alpha_k(\rho)^\frac{1}{2}\rho^2.
\end{equation}
Therefore for $\varepsilon\le\varepsilon_0$ we may apply the generalized Poincar\'e inequality \thref{Poincare} to the functions $f=\D_j u_k$ and $\delta=c\alpha_k(\rho)^\frac{1}{4}\rho$ to get a constant vector $\eta_k$ with $|\eta_k|\le c\varepsilon^\frac{1}{6}+\delta_k\le c$ such that
$$\int_{\Omega_k}|\D u_k-\eta_k|^2\le c\rho^2\int_{\Omega_k}|\D^2 u_k|^2+c\alpha_k(\rho)^\frac{1}{4}\rho^2\sup_{\Omega_k}|\D u_k|^2.$$
Since
$$\int_{\Omega_k}|\D^2 u_k|^2\le c\int_{\graph u_k}|\rm{A}_k|^2\,d\Hzwei\le c\int_{\Sigma_k\cap B_{2\rho}(\xi)}|\rm{A}_k|^2\,d\Hzwei\le c\alpha_k(\rho),$$
it follows for $\varepsilon\le\varepsilon_0$ that
\begin{equation}\label{22}
\int_{\Omega_k}\left|\D u_k-\eta_k\right|^2\le c\alpha_k(\rho)^\frac{1}{4}\rho^2.
\end{equation}
Let again $\overline{u_k}\in C^{1,1}(B_{\lambda_k}(\xi)\cap L, L^\perp)$ be an extension of $u_k$ to the hole disc $B_{\lambda_k}(\xi)\cap L$ as in \thref{extension}, i.e. $\overline{u_k}=u_k$ in $\Omega_k$. We again have that
\begin{eqnarray*}
\left\|\overline{u_k}\right\|_{L^\infty(B_{\lambda_k}(\xi)\cap L)} & \le & c\varepsilon^\frac{1}{6}\rho+\delta_k\hspace{0,2cm}\le\hspace{0,2cm}c,\\
\left\|\D\overline{u_k}\right\|_{L^\infty(B_{\lambda_k}(\xi)\cap L)} & \le & c\varepsilon^\frac{1}{6}+\delta_k\hspace{0,2cm}\le\hspace{0,2cm}c.
\end{eqnarray*}
From the gradient estimates for the function $\overline{u_k}$, since $|\eta_k|\le c$, from (\ref{21}), (\ref{22}) and the choice of $\rho_0$ we get that
\begin{eqnarray*}
\int_{B_{\lambda_k}(\xi)\cap L}\left|\D\overline{u_k}-\eta_k\right|^2 & = & \int_{\Omega_k}\left|\D u_k-\eta_k\right|^2+\sum_m\int_{d_{k,m}}\left|\D\overline{u_k}-\eta_k\right|^2 \\
& \le & c\alpha_k(\rho)^\frac{1}{4}\rho^2\phantom{\sum_m\int_{d_m^k}\left|\D\overline{u_k}-\eta_k\right|^2} \\
& \le & c\varepsilon^\frac{1}{2}\rho^2,
\end{eqnarray*}
and therefore in view of (\ref{20}) (we will always write $\alpha$ even if it might change from line to line) that
\begin{equation}\label{23}
\liminf_{k\to\infty}\int_{B_{\lambda_k}(\xi)\cap L}|\D\overline{u_k}-\eta_k|^2\le\min\left\{c\rho^{2+\alpha},c\varepsilon^\frac{1}{2}\rho^2\right\}\quad\text{for all }\rho\le\theta\frac{\rho_0}{64}.
\end{equation}
Since $\lambda_k>\frac{\rho}{4}$ the sequence $\overline{u_k}$ is therefore equicontinous and uniformly bounded in $C^1(B_\frac{\rho}{4}(\xi)\cap L,L^\perp)$ and $W^{1,2}(B_\frac{\rho}{4}(\xi)\cap L,L^\perp)$ and we get the existence of a function $u_\xi\in C^{0,1}(B_\frac{\rho}{4}(\xi)\cap L,L^\perp)$ such that (after passing to a subsequence)
\begin{eqnarray*}
\overline{u_k} & \to & u_\xi\quad\text{in }C^0(B_\frac{\rho}{4}(\xi)\cap L,L^\perp),\phantom{\frac{1}{\rho}} \\
\overline{u_k} & \rightharpoonup & u_\xi\quad\text{weakly in }W^{1,2}(B_\frac{\rho}{4}(\xi)\cap L,L^\perp),\phantom{\frac{1}{\rho}} \\
& & \hspace{-1,4cm}\frac{1}{\rho}\|u_\xi\|_{L^\infty(B_\frac{\rho}{4}(\xi)\cap L)}+\|\D u_\xi\|_{L^\infty(B_\frac{\rho}{4}(\xi)\cap L)} \le c\varepsilon^\frac{1}{6}.
\end{eqnarray*}
Remark: Be aware that the limit function depends on the point $\xi$, since our sequence comes (more or less) from the graphical decomposition lemma (which is a local statement) and therefore depends on $\xi$.\\[0,3cm]
Moreover we have that $\eta_k\to\eta$ with $|\eta|\le c\varepsilon^\frac{1}{6}$. Since $\D\overline{u_k}\rightharpoonup\D u_\xi$ weakly in $L^2(B_\frac{\rho}{4}(\xi)\cap L)$, it follows from lower semicontinuity and (\ref{23}) that
\begin{equation}\label{24}
\int_{B_\frac{\rho}{4}(\xi)\cap L}|\D u_\xi-\eta|^2\le c\rho^{2+\alpha}\le c\varepsilon^\frac{1}{2}\rho^2\quad\text{for all }\rho\le\theta\frac{\rho_0}{64}.
\end{equation}
In the next lemma we show that our limit varifold is given by a graph around the good points.

\begin{Le}\label{mu=graph}
For all $\xi\in\Sigma\cap B_\frac{\rho_0}{2}(\xi_0)$ and all $\rho<\theta\frac{\rho_0}{512}$ we have that
$$\mu\llcorner B_\rho(\xi)=\Hzwei\llcorner\left(\graph u_\xi\cap B_\rho(\xi)\right),$$
where $u_\xi\in C^{0,1}(B_\frac{\rho}{4}(\xi)\cap L,L^\perp)$ is as above.
\end{Le}

\begin{Pro}
From the definition of $\overline{u_k}$ it follows for $\rho\le\theta\frac{\rho_0}{64}$ that
\begin{eqnarray}\label{800}
\Hzwei\llcorner\left(\Sigma_k\cap B_\rho(\xi)\right) & = & \Hzwei\llcorner\left(D_k\cap B_\rho(\xi)\right)\nonumber \\
& = & \Hzwei\llcorner\left(\graph\overline{u_k}\cap B_\rho(\xi)\right)+\Hzwei\llcorner\left(D_k\setminus\graph\overline{u_k}\cap B_\rho(\xi)\right)\nonumber \\
& & -\Hzwei\llcorner\left(\graph\overline{u_k}\setminus D_k\cap B_\rho(\xi)\right)\nonumber \\
& = & \Hzwei\llcorner\left(\graph\overline{u_k}\cap B_\rho(\xi)\right)+\theta_k^\xi,\phantom{\left(\rho\le\frac{\rho_0}{32}\right)}
\end{eqnarray}
where $\theta_k^\xi$ is given by
$$\theta_k^\xi=\Hzwei\llcorner\left(D_k\setminus\graph\overline{u_k}\cap B_\rho(\xi)\right)-\Hzwei\llcorner\left(\graph\overline{u_k}\setminus D_k\cap B_\rho(\xi)\right)=\theta_k^1-\theta_k^2$$
is a signed measure. The total mass $|\theta_k^\xi|$ of $\theta_k^\xi$, namely $\theta_k^1(\Rdrei)+\theta_k^2(\Rdrei)$, can be estimated in view of (\ref{20}), (\ref{21}) and (\ref{21'}) by
\begin{eqnarray*}
\theta_k^1(\Rdrei)+\theta_k^2(\Rdrei) & \le & \sum_n\Hzwei\left(P_n^k\right)+\sum_m\int_{d_{k,m}}\sqrt{1+|\D\overline{u_k}|^2} \\
& \le & c\alpha_k(\rho)^\frac{1}{2}\rho^2\phantom{\sum_m\int_{d_m^k}\sqrt{1+|\D\overline{u_k}|^2}} \\
& \le & c\varepsilon\rho^2.\phantom{\sum_m\int_{d_m^k}\sqrt{1+|\D\overline{u_k}|^2}}
\end{eqnarray*}
It follows from (\ref{20}) that 
\begin{equation}\label{25}
\liminf_{k\to\infty}\left(\theta_k^1(\Rdrei)+\theta_k^2(\Rdrei)\right)\le c\rho^{2+\alpha}\le c\varepsilon\rho^2.
\end{equation}
By taking limits in the measure theoretic sense we get that
\begin{equation}\label{26}
\mu\llcorner B_\rho(\xi)=\Hzwei\llcorner\left(\graph u_\xi\cap B_\rho(\xi)\right)+\theta_\xi,
\end{equation}
where $\theta_\xi$ is a signed measure with total mass $|\theta_\xi|\le c\rho^{2+\alpha}\le c\varepsilon^\frac{1}{4}\rho^2$. This equation holds for all $\rho\le\theta\frac{\rho_0}{64}$ such that 
$$\mu\left(\partial B_\rho(\xi)\right)=\Hzwei\llcorner\graph u_\xi\left(\partial B_\rho(\xi)\right)=0,$$
which holds for a.e. $\rho\le\theta\frac{\rho_0}{64}$.\\[0,3cm]
To prove (\ref{26}) let $U\subset\Rdrei$ open.\\[0,3cm]
1.) Let $\rho\le\theta\frac{\rho_0}{64}$ such that $\mu\left(\partial B_\rho(\xi)\right)=0$. Moreover assume that $\mu\llcorner B_\rho(\xi)\left(\partial U\right)=0$. Therefore $\mu\left(\partial\left(U\cap B_\rho(\xi)\right)\right)=0$ and we get $\mu_k\left(U\cap B_\rho(\xi)\right)\to\mu\left(U\cap B_\rho(\xi)\right)$. It follows that
\begin{equation}
\Hzwei\llcorner\left(\Sigma_k\cap B_\rho(\xi)\right)(U)\to\mu\llcorner B_\rho(\xi)(U).
\end{equation}
2.) Let $\rho\le\theta\frac{\rho_0}{64}$ such that $\Hzwei\llcorner\graph u_\xi\left(\partial B_\rho(\xi)\right)=0$. Moreover assume that $\Hzwei\llcorner\left(\graph u_\xi\cap B_\rho(\xi)\right)\left(\partial U\right)=0$. We have that
$$\Hzwei\llcorner\left(\graph\overline{u_k}\cap B_\rho(\xi)\right)(U)=\int_{L}\chiup_{\phantom{}_{U\cap B_\rho(\xi)}}(x+\overline{u_k}(x))\sqrt{1+|\D\overline{u_k}(x)|^2}.$$
It follows from the $L^\infty$-bounds for $\overline{u_k}$ that
\begin{eqnarray*}
\left|\int_{L}\chiup_{\phantom{}_{U\cap B_\rho(\xi)}}(x+\overline{u_k}(x))\sqrt{1+|\D\overline{u_k}(x)|^2}-\int_{L}\chiup_{\phantom{}_{U\cap B_\rho(\xi)}}(x+u_\xi(x))\sqrt{1+|\D u_\xi(x)|^2}\right| & & \\
& \hspace{-22cm}\le & \hspace{-11cm}c\int_L\left|\chiup_{\phantom{}_{U\cap B_\rho(\xi)}}(x+\overline{u_k}(x))-\chiup_{\phantom{}_{U\cap B_\rho(\xi)}}(x+u_\xi(x))\right| \\
& & \hspace{-10cm}+\int_L\chiup_{\phantom{}_{U\cap B_\rho(\xi)}}(x+u_\xi(x))\left|\sqrt{1+|\D\overline{u_k}(x)|^2}-\sqrt{1+|\D u_\xi(x)|^2}\right|.
\end{eqnarray*}
Since $\overline{u_k}\to u_\xi$ uniformly and $\Hzwei\llcorner\graph u_\xi\left(\partial\left(U\cap B_\rho(\xi)\right)\right)=0$ we first of all get that
\begin{equation*}\label{27}
\chiup_{\phantom{}_{U\cap B_\rho(\xi)}}(x+\overline{u_k}(x))\to\chiup_{\phantom{}_{U\cap B_\rho(\xi)}}(x+u_\xi(x))\quad\text{for a.e. }x\in L.
\end{equation*}
The dominated convergence theorem yields
\begin{equation*}
\int_L\left|\chiup_{\phantom{}_{U\cap B_\rho(\xi)}}(x+\overline{u_k}(x))-\chiup_{\phantom{}_{U\cap B_\rho(\xi)}}(x+u_\xi(x))\right|\to0.
\end{equation*}
On the other hand we have that
\begin{eqnarray*}
\int_L\chiup_{\phantom{}_{U\cap B_\rho(\xi)}}(x+u_\xi(x))\left|\sqrt{1+|\D\overline{u_k}(x)|^2}-\sqrt{1+|\D u_\xi(x)|^2}\right| \phantom{\sqrt{|\D u(x)|}}& & \\
& \hspace{-18cm}\le & \hspace{-9cm}c\int_L\chiup_{\phantom{}_{U\cap B_\rho(\xi)}}(x+u_\xi(x))\left|\D\overline{u_k}(x)-\eta_k\right|+c\int_L\chiup_{\phantom{}_{U\cap B_\rho(\xi)}}(x+u_\xi(x))\left|\eta_k-\eta\right| \\
& & \hspace{-8cm}+c\int_L\chiup_{\phantom{}_{U\cap B_\rho(\xi)}}(x+u_\xi(x))\left|\eta-\D u_\xi(x)\right|.
\end{eqnarray*}
From the $L^\infty$-bound for $u_\xi$ it follows that $\chiup_{\phantom{}_{U\cap B_\rho(\xi)}}(x+u_\xi(x))=0$ if $x\notin B_{\left(1-c\varepsilon^\frac{1}{6}\right)\rho}(\xi)\cap L$ and we get that
$$\left(\int_L\chiup_{\phantom{}_{U\cap B_\rho(\xi)}}(x+u_\xi(x))\right)^\frac{1}{2}\le\Lzwei\left(B_{\left(1-c\varepsilon^\frac{1}{6}\right)\rho}(\xi)\cap L\right)^\frac{1}{2}\le c\rho.$$
In view of (\ref{23}), (\ref{24}) and since $\eta_k\to\eta$ we get that
$$\liminf_{k\to\infty}\int_L\chiup_{\phantom{}_{U\cap B_\rho(\xi)}}(x+u_\xi(x))\left|\sqrt{1+|\D\overline{u_k}(x)|^2}-\sqrt{1+|\D u_\xi(x)|^2}\right|\le c\rho^{2+\alpha}\le c\varepsilon^\frac{1}{4}\rho^2,$$
and it follows after all that
\begin{equation*}
\Hzwei\llcorner\left(\graph\overline{u_k}\cap B_\rho(\xi)\right)(U)=\Hzwei\llcorner\left(\graph u_\xi\cap B_\rho(\xi)\right)(U)+\tilde\theta_k^\xi(U),
\end{equation*}
where $\tilde\theta_k^\xi$ is a signed measure with $\liminf_{k\to\infty}|\tilde\theta_k^\xi|\le c\rho^{2+\alpha}\le c\varepsilon^\frac{1}{4}\rho^2$. After passing to a subsequence, the $\tilde\theta_k^\xi$'s converge to some signed measure $\tilde\theta_\xi$ with total mass $|\tilde\theta_\xi|\le c\rho^{2+\alpha}\le c\varepsilon^\frac{1}{4}\rho^2$. Assume that $\tilde\theta_\xi(\partial U)=0$. Then it follows that $\tilde\theta_k^\xi(U)\to\tilde\theta_\xi(U)$ and therefore we get
\begin{equation}
\lim_{k\to\infty}\Hzwei\llcorner\left(\graph\overline{u_k}\cap B_\rho(\xi)\right)(U)=\Hzwei\llcorner\left(\graph u_\xi\cap B_\rho(\xi)\right)(U)+\tilde\theta_\xi(U).
\end{equation}
3.) Since the $\theta_k^\xi$'s were signed measures such that $\liminf|\theta_k^\xi|\le c\rho^{2+\alpha}\le c\varepsilon^\frac{1}{4}\rho^2$, they converge in the weak sense (after passing to a subsequence) to a signed measure $\overline{\theta}_\xi$ with total mass $|\overline{\theta}_\xi|\le c\rho^{2+\alpha}\le c\varepsilon^\frac{1}{4}\rho^2$. Assuming that $\overline{\theta}_\xi(\partial U)=0$, we get $\theta_k^\xi(U)\to\overline{\theta}_\xi(U)$. \\[0,3cm]
Now by taking limits in (\ref{800}) it follows that
\begin{equation}
\mu\llcorner B_\rho(\xi)(U)=\Hzwei\llcorner\left(\graph u_\xi\cap B_\rho(\xi)\right)(U)+\theta_\xi(U),
\end{equation}
where $\theta_\xi=\overline{\theta}_\xi+\tilde\theta_\xi$ is a signed measure with total mass $|\theta_\xi|\le c\rho^{2+\alpha}\le c\varepsilon^\frac{1}{4}\rho^2$. Notice that this equation holds for every $U\subset\Rdrei$ open such that
$$\mu\llcorner B_\rho(\xi)(\partial U)=\Hzwei\llcorner\left(\graph u_\xi\cap B_\rho(\xi)\right)(\partial U)=\overline\theta_\xi(\partial U)=\tilde\theta_\xi(\partial U)=0.$$
By choosing an appropriate exhaustion this equation holds for arbitrary open sets $U\subset\Rdrei$ and (\ref{26}) follows.\\[0,3cm]
Now choose a radius $\rho\in\left(\theta\frac{\rho_0}{128},\theta\frac{\rho_0}{64}\right)$ such that (\ref{26}) holds. We take a closer look to two cases.\\[0,3cm]
1.) Let $x\in\Sigma\cap B_\frac{\rho}{2}(\xi)$: Notice that by (\ref{4'}) and the choice of $\rho_0$
$$\Will(\mu,B_\frac{\rho}{2}(x))\le\liminf_{k\to\infty}\Will(\mu_k,B_\frac{\rho}{2}(x))\le c\int_{\Sigma_k\cap B_\rho(\xi)}|\rm{A}_k|^2\,d\Hzwei\le 2\varepsilon^2.$$
Since $\theta^2(\mu,\cdot)\ge1$ on $\spt\mu$, it follows for $\varepsilon\le\varepsilon_0$ from \thref{Monotonie} that
$$\mu\llcorner B_\rho(\xi)\left(B_\frac{\rho}{2}(x)\right)=\mu\left(B_\frac{\rho}{2}(x)\right)\ge c\rho^2.$$
From (\ref{26}), especially the bound on the total mass of $\theta_\xi$, it follows that
$$c\rho^2\le\Hzwei\left(\graph u_\xi\cap B_\frac{\rho}{2}(x)\right)+c\varepsilon^\frac{1}{4}\rho^2.$$
Therefore $\Hzwei\left(\graph u_\xi\cap B_\frac{\rho}{2}(x)\right)>0$ for $\varepsilon\le\varepsilon_0$ and thus $x\in\graph u_\xi$.\\[0,3cm]
2.) Let $x\in\graph u_\xi\cap B_\frac{\rho}{2}(\xi)$: Write $x=z+u_\xi(z)$. If $y\in B_\frac{\rho}{4}(z)\cap L$, it follows from the estimates for $u_\xi$ that $y+u_\xi(y)\in B_\frac{\rho}{2}(x)$ for $\varepsilon\le\varepsilon_0$. Therefore we get that
$$\Hzwei\llcorner\graph u_\xi\left(B_\frac{\rho}{2}(x)\right)\ge\int\chiup_{\phantom{}_{B_\frac{\rho}{2}(x)}}(y+u_\xi(y))\ge c\rho^2.$$
As above it follows that $x\in\Sigma$ for $\varepsilon\le\varepsilon_0$.\\[0,3cm]
After all we get for $\varepsilon\le\varepsilon_0$ 
\begin{equation}\label{37}
\Sigma\cap B_\rho(\xi)=\graph u_\xi\cap B_\rho(\xi)\quad\text{for all }\rho<\theta\frac{\rho_0}{256}.
\end{equation}
Moreover we get that the function $u_\xi$ does not depend on the point $\xi$ in the following sense: Let $x\in\Sigma\cap B_\frac{\rho_0}{2}(\xi_0)$ and $\tau<\theta\frac{\rho_0}{256}$. Then we have that 
\begin{equation}\label{38}
\graph u_\xi\cap\left(B_\rho(\xi)\cap B_\tau(x)\right)=\graph u_x\cap\left(B_\rho(\xi)\cap B_\tau(x)\right).
\end{equation}
In the next step choose $\rho\in\left(\theta\frac{\rho_0}{512},\theta\frac{\rho_0}{256}\right)$ such that
$$\mu\left(\partial B_\rho(\xi)\right)=\Hzwei\llcorner\graph u_\xi\left(\partial B_\rho(\xi)\right)=0,$$
and that therefore due to (\ref{26})
\begin{equation}\label{39}
\mu\llcorner B_\rho(\xi)=\Hzwei\llcorner\left(\graph u_\xi\cap B_\rho(\xi)\right)+\theta_\xi.
\end{equation}
Let $x\in\Sigma\cap B_\rho(\xi)=\graph u_\xi\cap B_\rho(\xi)$ and $\tau>0$ such that $B_\tau(x)\subset B_\rho(\xi)$ and such that (due to (\ref{26}) for the point $x$)
\begin{equation}\label{40}
\mu\llcorner B_\tau(x)=\Hzwei\llcorner\left(\graph u_x\cap B_\tau(x)\right)+\theta_x,
\end{equation}
where $\theta_x$ is a signed measure with total mass smaller than $c\tau^{2+\alpha}$.\\[0,3cm]
From (\ref{38}), (\ref{39}) and (\ref{40}) it follows that 
$$\theta_\xi\left(B_\tau(x)\right)=\theta_x\left(B_\tau(x)\right)$$
and we get a nice decay for the signed measure $\theta_\xi$, namely
\begin{equation}\label{41}
\lim_{\tau\to0}\frac{\theta_\xi\left(B_\tau(x)\right)}{\tau^2}=0.
\end{equation}
Since we already know that $\theta^2(\mu,\cdot)\ge1$ on $\Sigma$, it follows from (\ref{39}) that
$$\frac{\Hzwei\llcorner\left(\graph u_\xi\cap B_\rho(\xi)\right)\left(B_\tau(x)\right)}{\mu\llcorner B_\rho(\xi)\left(B_\tau(x)\right)}=1-\frac{\theta_\xi\left(B_\tau(x)\right)}{\mu\llcorner B_\rho(\xi)\left(B_\tau(x)\right)},$$
and by (\ref{41}) the right hand side goes to 1 for $\tau\to0$. This shows that
$$\D_{\mu\llcorner B_\rho(\xi)}\left(\Hzwei\llcorner\left(\graph u_\xi\cap B_\rho(\xi)\right)\right)(x)=1$$ 
for all $x\in\Sigma\cap B_\rho(\xi)=\graph u_\xi\cap B_\rho(\xi)$ and the lemma follows from the theorem of Radon-Nikodym.
\end{Pro}

Now let $\xi_0\in\Sigma\setminus\B$. Since we already know that $\mu$ admits a generalized mean curvature vector $\vec{H}\in L^2(\mu)$, it follows from the definition of the weak mean curvature vector and by applying \thref{mu=graph} to $\xi_0$ (and writing $u$ for $u_{\xi_0}$) that $u$ is a weak solution of the mean curvature equation, namely $u$ is a weak solution in $W^{1,2}_0\left(B_{\theta\frac{\rho_0}{512}}(\xi_0)\cap L,L^\perp\right)$ of
\begin{equation}\label{28}
\sum_{i,j=1}^2\partial_j\left(\sqrt{\det g}\,\,g^{ij}\partial_i F\right)=\sqrt{\det g}\,\,\vec{H}\circ F,
\end{equation}
where $F(x)=x+u(x)$ and $g_{ij}=\delta_{ij}+\partial_i u\cdot\partial_j u$.\\[0,3cm]
Since the norm of the mean curvature vector can be estimated by the norm of the second fundamental form, it follows from \thref{2ff-absch} and (\ref{4'}) applied to $B_\rho(\xi)$ that for all $\xi\in B_\frac{\rho_0}{2}(\xi_0)$ and all $\rho\le\theta\frac{\rho_0}{128}$
$$\int_{B_\rho(\xi)}|\vec{H}|^2\,d\mu\le c\rho^\alpha.$$
\thref{mu=graph} yields $\mu\llcorner B_\rho(\xi)=\Hzwei\llcorner\left(\graph u\cap B_\rho(\xi)\right)$ for all points $\xi\in B_{\theta\frac{\rho_0}{1024}}(\xi_0)$ and all $\rho\le\theta\frac{\rho_0}{1024}$ and therefore 
\begin{equation}\label{29}
\int_{\graph u\cap B_\rho(\xi)}|\vec{H}|^2\,d\Hzwei\le c\rho^\alpha
\end{equation}
for all $\xi\in B_{\theta\frac{\rho_0}{1024}}(\xi_0)$ and all $\rho\le\theta\frac{\rho_0}{1024}$.\\[0,3cm]
Using a standard difference quotient argument (as for example in \cite{GT}, Theorem 8.8), it follows from (\ref{28}) and $\vec{H}\in L^2(\mu)$ that 
\begin{equation*}
u\in W^{2,2}_{loc}\left(B_{\theta\frac{\rho_0}{1024}}(\xi_0)\cap L,L^\perp\right).
\end{equation*}
Now let $\varphi=\eta\nu$, where $\nu\perp L$ and $\eta\in W^{1,2}_0\left(B_{\theta\frac{\rho_0}{1024}}(\xi_0)\cap L\right)$ is of the form 
$$\eta=f^2\left(\partial_l u-\fint_{B_\rho(x)\setminus B_\frac{\rho}{2}(x)\cap L}\partial_l u\right),$$
where $x\in B_{\theta\frac{\rho_0}{2048}}(\xi_0)\cap L$, $\rho\le\theta\frac{\rho_0}{2048}$ and $f\in C_c^\infty(B_\rho(x)\cap L)$ such that $0\le f\le1$, $f\equiv1$ on $B_{\frac{\rho}{2}}(x)\cap L$ and $|\D f|\le\frac{c}{\rho}$.\\[0,3cm]
By applying (\ref{28}) to this test functions $\varphi$ we get in view of (\ref{29}) that
\begin{equation}\label{30}
\int_{B_\rho(x)\cap L}|\D^2 u|^2\le c\rho^\alpha\quad\text{for all }x\in B_{\theta\frac{\rho_0}{2048}}(\xi_0)\cap L\text{ and all }\rho\le\theta\frac{\rho_0}{2048}.
\end{equation}
From Morrey's lemma (see \cite{GT}, Theorem 7.19) it follows that
\begin{equation}\label{31}
u\in C^{1,\alpha}\left(B_{\theta\frac{\rho_0}{2048}}(\xi_0)\cap L,L^\perp\right)\cap W^{2,2}\left(B_{\theta\frac{\rho_0}{2048}}(\xi_0)\cap L,L^\perp\right).
\end{equation}
Thus we have shown that our limit varifold $\Sigma$ can be written as a $C^{1,\alpha}\cap W^{2,2}$-graph away from the bad points. \\[0,5cm]
Now we will handle the bad points $\B$ and prove a similar power decay as in \thref{2ff-absch} for balls around the bad points. Since the bad points are discrete and since we want to prove a local decay, we assume that there is only one bad point $\xi$.\\[0,3cm]
As mentioned in the definition of the bad points (see (\ref{7})), the radon measures $\alpha_k=\mu_k\llcorner|\rm{A}_k|^2$ converge weakly to a radon measure $\alpha$, and it follows for all $z\in\Rdrei$ that $\alpha(B_\rho(z)\setminus\{z\})\to0$ for $\rho\to0.$ Therefore for given $\varepsilon>0$ there exists a $\rho_0>0$ such that
$$\alpha(B_\rho(\xi)\setminus\{\xi\})<\varepsilon^2\quad\text{for all }\rho\le\rho_0.$$
Since $\alpha_k\to\alpha\text{ in }C_c^0(\Rdrei)'$, it follows that for $\rho<\rho_0$ and $k\in\varmathbb{N}$ sufficiently large
\begin{equation}\label{32}
\int_{\Sigma_k\cap B_\rho(\xi)\setminus B_\frac{\rho}{2}(\xi)}|\rm{A}_k|^2\,d\Hzwei<\varepsilon^2.
\end{equation}

Moreover it follows from \thref{Monotonie} applied to our minimizing sequence $\Sigma_k$ and (\ref{2}) that for all $\sigma>0$
$$\int_{\Rdrei\setminus B_\sigma(\xi)}\left|\frac{1}{4}\vec{H}_k(x)+\frac{(x-\xi)^\perp}{|x-\xi|^2}\right|^2\,d\mu_k(x)\le\frac{1}{4\pi}\Will(\Sigma_k)-\theta^2(\mu_k,\xi)\le c,$$
where $c$ is an universal constant independent of $k$ and $\sigma$. Here $\perp$ denotes the projection onto $T_x\Sigma_k$. Rewriting the left hand side and using Cauchy-Schwarz we get
$$\int\chiup_{\Rdrei\setminus B_\sigma(\xi)}\frac{|(x-\xi)^\perp|^2}{|x-\xi|^4}\,d\mu_k(x)\le c,$$
where again $c$ is an universal constant independent of $k$ and $\sigma$. Now we can use the monotone convergence theorem to get for $\sigma\to0$ that the integral
$$\int\frac{|(x-\xi)^\perp|^2}{|x-\xi|^4}\,d\mu_k(x)$$
exists for all $k$ and is bounded by a uniform, universal constant $c$ independent of $k$. Moreover the function 
$$f_k=\frac{|(x-\xi)^\perp|^2}{|x-\xi|^4}\in L^1(\mu_k).$$
Now define the radon measures
$$\beta_k=f_k\llcorner\mu_k.$$
It follows that $\beta_k(\Rdrei)\le c$ and therefore (after passing to a subsequence) there exists a radon measure $\beta$ such that $\beta_k\to\beta$ in $C_c^0(\Rdrei)'$. Moreover $\beta(B_\rho(\xi)\setminus\{\xi\})\to0$ for $\rho\to0$. Therefore there exists a $\rho_0$ such that 
$$\beta(B_\rho(\xi)\setminus\{\xi\})<\varepsilon^4\quad\text{for all }\rho\le\rho_0.$$
Let $\rho<\rho_0$ and $g\in C_c^0(B_{\rho_0}(\xi)\setminus\{\xi\})$ such that $0\le g\le1$ and $g\ge\chiup_{B_\rho(\xi)\setminus B_\frac{\rho}{2}(\xi)}$. It follows that
$$\int\chiup_{B_\rho(\xi)\setminus B_\frac{\rho}{2}(\xi)}\frac{|(x-\xi)^\perp|^2}{|x-\xi|^4}\,d\mu_k(x)\le\int g\,d\beta_k\to\int g\,d\beta\le\beta(B_{\rho_0}(\xi)\setminus\{\xi\})<\varepsilon^4.$$
Thus we get for $k$ sufficiently large that
$$\int_{\Sigma_k\cap B_\rho(\xi)\setminus B_\frac{\rho}{2}(\xi)}\frac{|(x-\xi)^\perp|^2}{|x-\xi|^4}\,d\Hzwei(x)\le\varepsilon^4.$$
Now let $\mathcal{B}_k=\left\{x\in\Sigma_k\cap B_{\rho_0}(\xi)\,\Big|\,\frac{|(x-\xi)^\perp|}{|x-\xi|}>\varepsilon\right\}$. It follows for $\rho<\rho_0$ and $k\in\varmathbb{N}$ sufficiently large that
\begin{equation}\label{33}
\varepsilon\Hzwei\left(\Sigma_k\cap B_\rho(\xi)\setminus B_\frac{\rho}{2}(\xi)\cap\mathcal{B}_k\right)\le c\varepsilon^2\rho^2.
\end{equation}
Moreover by choosing $\rho_0\le\frac{2}{3\sqrt{8\pi}}$ we also get for $\rho<\rho_0$ and for $k$ large that
\begin{equation}\label{34}
\Sigma_k\cap\partial B_{\frac{3}{4}\rho}(\xi)\neq\emptyset.
\end{equation}
To prove this notice that due to the diameter estimate in Lemma 1.1 in \cite{SL} we have
$$\diam\Sigma_k\ge\sqrt{\frac{\Hzwei(\Sigma_k)}{\Will(\Sigma_k)}}\ge\sqrt{\frac{1}{8\pi}}.$$
Let $\xi_k\in\Sigma_k$ such that $\xi_k\to\xi$. It follows that $\Sigma_k\cap B_{\frac{3}{4}\rho}(\xi)\neq\emptyset$ for $k\in\varmathbb{N}$ sufficiently large. Now suppose that $\Sigma_k\cap\partial B_{\frac{3}{4}\rho}(\xi)=\emptyset$. Since $\Sigma_k$ is connected, we get that $\Sigma_k\subset B_{\frac{3}{4}\rho}(\xi)$ and therefore $\diam\Sigma_k\le\frac{3}{2}\rho<\frac{3}{2}\rho_0\le\frac{1}{\sqrt{8\pi}}$, a contradiction.\\[0,3cm]
After all according to (\ref{32})-(\ref{34}) the following is shown: For $\rho<\rho_0$ and $k\in\varmathbb{N}$ sufficiently large we have that
\begin{eqnarray*}
& (i) & \int_{\Sigma_k\cap B_\rho(\xi)\setminus B_\frac{\rho}{2}(\xi)}|\rm{A}_k|^2\,d\Hzwei<\varepsilon^2,\\
& (ii) & \frac{|(x-\xi)^\perp|}{|x-\xi|}\le\varepsilon\quad\text{for all }x\in\left(\Sigma_k\cap B_\rho(\xi)\setminus B_\frac{\rho}{2}(\xi)\right)\setminus\mathcal{B}_k,\\
& & \text{where }\mathcal{B}_k\subset\Sigma_k\cap B_{\rho_0}(\xi)\text{ with }\Hzwei\left(\Sigma_k\cap B_\rho(\xi)\setminus B_\frac{\rho}{2}(\xi)\cap\mathcal{B}_k\right)\le c\varepsilon\rho^2\\
& & \text{and $(x-\xi)^\perp=(x-\xi)-P_{T_x\Sigma_k}(x-\xi)$,}\\
& (iii) & \Sigma_k\cap\partial B_{\frac{3}{4}\rho}(\xi)\neq\emptyset.\phantom{\int_{\Sigma_k\cap B_\rho(\xi)\setminus B_\frac{\rho}{2}(\xi)}}
\end{eqnarray*}
Let $z_k\in\Sigma_k\cap\partial B_{\frac{3}{4}\rho}(\xi)$. It follows that
$$\int_{\Sigma_k\cap B_\frac{\rho}{8}(z_k)}|\rm{A}_k|^2\,d\Hzwei\le\int_{\Sigma_k\cap B_\rho(\xi)\setminus B_\frac{\rho}{2}(\xi)}|\rm{A}_k|^2\,d\Hzwei\le\varepsilon^2.$$
The Monotonicity formula applied to $z_k$ and $\Sigma_k$ yields that we may apply the graphical decomposition lemma to $\Sigma_k$, $z_k$ and infinitely many $k\in\varmathbb{N}$ as well as Lemma 1.4 in \cite{SL} to get as in \thref{final} that there exists a $\theta\in\left(0,\frac{1}{2}\right)$ (independent of $j\in\{1,\ldots,P\}$ and $k\in\varmathbb{N}$) and pairwise disjoint subsets $P_1^k,\ldots,P_{N_k}^k\subset\Sigma_k$ such that 
$$\Sigma_k\cap\overline{B_{\theta\frac{\rho}{32}}(z_k)}=\left(\graph u_k\cup\bigcup_n P_n^k\right)\cap\overline{B_{\theta\frac{\rho}{32}}(z_k)},$$
where the following holds:  
\begin{enumerate}
\item The sets $P_n^k$ are closed topological discs disjoint from $\graph u_k$.
\item $u_k\in C^\infty(\overline{\Omega_k},L_k^\perp)$, where $L_k\subset\Rdrei$ is a 2-dim. plane such that $z_k\in L_k$ and $\Omega_k=\left(B_{\lambda_k}(z_k)\cap L_k\right)\setminus\bigcup_m d_{k,m}$, where $\lambda_k>\frac{\rho}{16}$ and where the sets $d_{k,m}$ are pairwise disjoint closed discs in $L_k$.
\item The following inequalities hold:
\begin{eqnarray}
& & \hspace{-1cm}\sum_m\diam d_{k,m}+\sum_n\diam P_n^k\le c\varepsilon^\frac{1}{2}\rho,\phantom{\sum_n\frac{1}{\rho}}\\
& & \hspace{-1cm}\|u_k\|_{L^\infty(\Omega_k)}\le c\varepsilon^{\frac{1}{6}}\rho+\delta_k\quad\text{where }\lim_{k\to\infty}\delta_k=0,\phantom{\frac{1}{\rho}\sum_n} \\
& & \hspace{-1cm}\|\D u_k\|_{L^\infty(\Omega_k)}\le c\varepsilon^{\frac{1}{6}}+\delta_k\quad\text{where }\lim_{k\to\infty}\delta_k=0.\phantom{\frac{1}{\rho}\sum_n}
\end{eqnarray}
\end{enumerate} 

In the next step we show that
\begin{equation}\label{35}
\dist\left(\xi,L_k\right)\le c\left(\varepsilon^\frac{1}{6}+\delta_k\right)\rho.
\end{equation}
To prove this notice first of all that it follows from \thref{Monotonie} applied to $z_k$, $\Sigma_k$ and (i) above that for $\varepsilon\le\varepsilon_0$
\begin{equation}\label{36}
\Hzwei(\Sigma_k\cap B_{\theta\frac{\rho}{32}}(z_k))\ge c\rho^2\quad\text{with $c$ independent of $k$}.
\end{equation}
Now to prove (\ref{35}) notice that
$$\left(\graph u_k\cap B_{\theta\frac{\rho}{32}}(z_k)\right)\setminus\mathcal{B}_k\neq\emptyset,$$
where $\mathcal{B}_k\subset\Sigma_k\cap B_{\rho_0}(\xi)$ is the set in (ii) above. This follows from the graphical decomposition above, the diameter estimates for the sets $P_n^k$, the area estimate concerning the set $\mathcal{B}_k$ in (ii) and (\ref{36}).\\[0,3cm]
Let $z\in\left(\graph u_k\cap B_{\theta\frac{\rho}{32}}(z_k)\right)\setminus\mathcal{B}_k\subset\left(\Sigma_k\cap B_\rho(\xi)\setminus B_\frac{\rho}{2}(\xi)\right)\setminus\mathcal{B}_k$. It follows from (ii) that
$$|\xi-\pi_{(z+T_z\Sigma_k)}(\xi)|\le\varepsilon|z-\xi|\le\varepsilon\left(|z-z_k|+|z_k-\xi|\right)\le c\varepsilon\rho.$$
Define the perturbed 2-dim. plane $\tilde L_k$ by $\tilde L_k=L_k+(z-\pi_{L_k}(z))$, where we have that $\dist(\tilde L_k,L_k)=|z-\pi_{L_k}(z)|\le c\varepsilon^\frac{1}{6}\rho$ (since $z\in\graph u_k\cap B_{\theta\frac{\rho}{32}}(z_k)$). Now it follows from Pythagoras that $|z-\pi_{\tilde L_k}(\pi_{(z+T_z\Sigma_k)}(\xi))|^2\le|z-\pi_{(z+T_z\Sigma_k)}(\xi)|^2\le|z-\xi|^2\le c\rho^2$. Since $z+T_z\Sigma_k$ can be parametrized in terms of $\D u_k(z)$ over $\tilde L_k$, we get that 
$$|\pi_{(z+T_z\Sigma_k)}(\xi)-\pi_{\tilde L_k}(\pi_{(z+T_z\Sigma_k)}(\xi))|\le\|\D u_k\|_{L^\infty}|z-\pi_{\tilde L_k}(\pi_{(z+T_z\Sigma_k)}(\xi))|\le c\left(\varepsilon^\frac{1}{6}+\delta_k\right)\rho.$$
Therefore we finally get that
\begin{eqnarray*}
\dist(\xi,L_k) & = & \left|\xi-\pi_{L_k}(\xi)\right| \\
& \le & \left|\xi-\pi_{L_k}(\pi_{(z+T_z\Sigma_k)}(\xi))\right|\\
& \le & \left|\xi-\pi_{(z+T_z\Sigma_k)}(\xi)\right|+\left|\pi_{(z+T_z\Sigma_k)}(\xi)-\pi_{\tilde L_k}(\pi_{(z+T_z\Sigma_k)}(\xi))\right| \\
& & +\left|\pi_{\tilde L_k}(\pi_{(z+T_z\Sigma_k)}(\xi))-\pi_{L_k}(\pi_{(z+T_z\Sigma_k)}(\xi))\right| \\
& \le & c\left(\varepsilon^\frac{1}{6}+\delta_k\right)\rho,
\end{eqnarray*}
and (\ref{35}) is shown.\\[0,3cm]
Since $\dist(\xi,L_k)\le c\left(\varepsilon^\frac{1}{6}+\delta_k\right)\rho$, we may assume (after translation) that $\xi\in L_k$ for all $k\in\varmathbb{N}$ and keeping the estimates for $u_k$. Moreover we again have that $L_k\to L=2\text{-dim. plane with }\xi\in L.$ Therefore for $k\in\varmathbb{N}$ sufficiently large we may assume that $L_k$ is a fixed 2-dim. plane $L$. \\[0,3cm]
Define the set
$$T_k=\left\{\tau\in\left(\theta\frac{\rho}{64},\theta\frac{\rho}{\sqrt{2}\cdot32}\right)\,\Bigg|\,\partial B_\tau(z_k)\cap\bigcup_m d_{k,m}=\emptyset\right\}.$$
It follows from the diameter estimates and the selection principle in \cite{SL} that for $\varepsilon\le\varepsilon_0$ there exists a $\tau\in\left(\theta\frac{\rho}{64},\theta\frac{\rho}{\sqrt{2}\cdot32}\right)$ such that $\tau\in T_k$ for infinitely many $k\in\varmathbb{N}$.\\[0,3cm]
Since $\xi\in L$ it follows from the choice of $\tau$ that for $\varepsilon\le\varepsilon_0$
$$\partial B_{\frac{3}{4}\rho}(\xi)\cap\partial B_\tau(z_k)\cap L=\left\{p_{1,k},p_{2,k}\right\},$$
where $p_{1,k}, p_{2,k}\in\left(B_{\theta\frac{\rho}{\sqrt{2}\cdot32}}(z_k)\cap L\right)\setminus\bigcup_m d_{k,m}$ are distinct.\\[0,3cm]
Define the image points $z_{i,k}\in\graph u_k$ by
$$z_{i,k}=p_{i,k}+u_k(p_{i,k}).$$
Using the $L^\infty$-estimates for $u_k$ we get for $\varepsilon\le\varepsilon_0$ that $\frac{5}{8}\rho<|z_{i,k}-\xi|<\frac{7}{8}\rho$ and therefore
$$\int_{\Sigma_k\cap B_\frac{\rho}{8}(z_{i,k})}|\rm{A}_k|^2\,d\Hzwei\le\int_{\Sigma_k\cap B_\rho(\xi)\setminus B_\frac{\rho}{2}(\xi)}|\rm{A}_k|^2\,d\Hzwei<\varepsilon^2.$$
Therefore we can again apply the graphical decomposition lemma to the points $z_{i,k}$. Thus there exist pairwise disjoint subsets $P_1^{i,k},\ldots,P_{N_{i,k}}^{i,k}\subset\Sigma_k$ such that 
$$\Sigma_k\cap\overline{B_{\theta\frac{\rho}{32}}(z_{i,k})}=\left(\graph u_{i,k}\cup\bigcup_n P_n^{i,k}\right)\cap\overline{B_{\theta\frac{\rho}{32}}(z_{i,k})},$$
where the following holds:  
\begin{enumerate}
\item The sets $P_n^{i,k}$ are closed topological discs disjoint from $\graph u_{i,k}$.
\item $u_{i,k}\in C^\infty(\overline{\Omega_{i,k}},L_{i,k}^\perp)$, where $L_{i,k}\subset\Rdrei$ is a 2-dim. plane such that $z_{i,k}\in L_{i,k}$ and $\Omega_{i,k}=\left(B_{\lambda_{i,k}}(z_{i,k})\cap L_{i,k}\right)\setminus\bigcup_m d_{i,k,m}$, where $\lambda_{i,k}>\frac{\rho}{16}$ and where the sets $d_{i,k,m}$ are pairwise disjoint closed discs in $L_{i,k}$.
\item The following inequalities hold:
\begin{eqnarray}
& & \hspace{-1cm}\sum_m\diam d_{i,k,m}+\sum_n\diam P_n^{i,k}\le c\varepsilon^\frac{1}{2}\rho,\phantom{\sum_n\frac{1}{\rho}}\\
& & \hspace{-1cm}\|u_{i,k}\|_{L^\infty(\Omega_{i,k})}\le c\varepsilon^{\frac{1}{6}}\rho+\delta_{i,k}\quad\text{where }\lim_{k\to\infty}\delta_{i,k}=0,\phantom{\sum_n\frac{1}{\rho}} \\
& & \hspace{-1cm}\|\D u_{i,k}\|_{L^\infty(\Omega_{i,k})}\le c\varepsilon^{\frac{1}{6}}+\delta_{i,k}\quad\text{where }\lim_{k\to\infty}\delta_{i,k}=0.\phantom{\sum_n\frac{1}{\rho}}
\end{eqnarray}
\vspace{0,3cm}
\end{enumerate} 
Since $\dist(z_{i,k},L)\le c\varepsilon^\frac{1}{6}\rho+\delta_k$ (this follows since $z_{i,k}\in\graph u_k$) and since the $L^\infty$-norms of $u_k$ and $u_{i,k}$ are small, we may assume (after translation and rotation as done before) that $L_{i,k}=L$.\\[0,3cm]
By continuing with this procedure we get after a finite number of steps, depending not on $\rho$ and $k\in\varmathbb{N}$, an open cover of $\partial B_{\frac{3}{4}\rho}(\xi)\cap L$ which also covers the set
$$B=\left\{x\in L\,\Big|\,\dist\left(x,\partial B_{\frac{3}{4}\rho}(\xi)\cap L\right)<\theta\frac{\rho}{\sqrt{2}\cdot64}\right\}$$
and which include finitely many, closed discs $d_{k,m}$ with
$$\sum_m\diam d_{k,m}\le c\varepsilon^\frac{1}{2}\rho.$$
We may assume that these discs are pairwise disjoint since otherwise we can exchange two intersecting discs by one disc whose diameter is smaller than the sum of the diameters of the intersecting discs.\\[0,3cm]
Because of the diameter estimate and again the selection principle there exists a $\tau\in\left(\theta\frac{\rho}{128},\theta\frac{\rho}{\sqrt{2}\cdot64}\right)$ such that
$$\left\{x\in L\,\Big|\,\dist\left(x,\partial B_{\frac{3}{4}\rho}(\xi)\cap L\right)=\tau\right\}\cap\bigcup_m d_{k,m}=\emptyset.$$
Finally we get the following: There exist pairwise disjoint subsets $P_1^k,\ldots,P_{N_k}^k\subset\Sigma_k$ such that 
$$\Sigma_k\cap\mathcal{A}(\rho)=\left(\graph u_k\cup\bigcup_n P_n^k\right)\cap\mathcal{A}(\rho),$$
where the following holds:  
\begin{enumerate}
\item The sets $P_n^k$ are closed topological discs disjoint from $\graph u_k$.
\item $u_k\in C^\infty(A_k(\rho),L^\perp)$, where $L\subset\Rdrei$ is a 2-dim. plane with $\xi\in L$.
\item The set $A_k(\rho)$ is given by
$$A_k(\rho)=\left\{x\in L\,\Big|\,\dist\left(x,\partial B_{\frac{3}{4}\rho}(\xi)\cap L\right)<\tau\right\}\setminus\bigcup_m d_{k,m},$$
where $\tau\in\left(\theta\frac{\rho}{128},\theta\frac{\rho}{\sqrt{2}\cdot64}\right)$ and where the sets $d_{k,m}$ are pairwise disjoint closed discs in $L$ which do not intersect $\left\{x\in L\,\Big|\,\dist\left(x,\partial B_{\frac{3}{4}\rho}(\xi)\cap L\right)=\tau\right\}$.
\item The set $\mathcal{A}(\rho)$ is given by
$$\mathcal{A}(\rho)=\left\{x+y\in\Rdrei\,\Big|\,x\in L, \dist\left(x,\partial B_{\frac{3}{4}\rho}(\xi)\cap L\right)<\tau, y\in L^\perp, |y|<\theta\frac{\rho}{64}\right\}.$$
\item The following inequalities hold:
\begin{eqnarray}
& & \hspace{-2,1cm}\sum_m\diam d_{k,m}+\sum_n\diam P_n^k\le c\varepsilon^\frac{1}{2}\rho,\phantom{\lim_{k\to\infty}\varepsilon^{\frac{1}{6}}}\\
\|u_k\|_{L^\infty(A_k(\rho))} & \le & c\varepsilon^{\frac{1}{6}}\rho+\delta_k\quad\text{where }\lim_{k\to\infty}\delta_k=0,\phantom{\sum_{j=1}\diam P_j^k} \\
\|\D u_k\|_{L^\infty(A_k(\rho))} & \le & c\varepsilon^{\frac{1}{6}}+\delta_k\quad\text{where }\lim_{k\to\infty}\delta_k=0.\phantom{\sum_{j=1}\diam P_j^k}
\end{eqnarray}
\end{enumerate}
From the estimates for the function $u_k$ and the diameter estimates for the sets $P_n^k$ we also get for $\varepsilon\le\varepsilon_0$ and $k$ sufficiently large that
\begin{equation*}
\Sigma_k\cap\mathcal{A}(\rho)\subset\left\{x+y\in\Rdrei\,\Big|\,x\in L, \dist\left(x,\partial B_{\frac{3}{4}\rho}(\xi)\cap L\right)<\tau, y\in L^\perp, |y|<\theta\frac{\rho}{128}\right\}.
\end{equation*}
Since $\Sigma_k\to\Sigma$ in the Hausdorff distance sense it follows that
$$\emptyset\neq\Sigma\cap\mathcal{A}(\rho)\subset\left\{x+y\in\Rdrei\Big|\,x\in L, \dist\left(x,\partial B_{\frac{3}{4}\rho}(\xi)\right)<\tau, y\in L_k^\perp, |y|<\theta\frac{\rho}{128}\right\}.$$
Now we show that for all $\rho<\rho_0$ (after choosing $\rho_0$ smaller if necessary)
$$\Sigma\cap\mathcal{A}(\rho)\cap B_{\left(\frac{3}{4}+\frac{\theta}{256}\right)\rho}(\xi)\setminus B_{\left(\frac{3}{4}-\frac{\theta}{256}\right)\rho}(\xi)=\Sigma\cap B_{\left(\frac{3}{4}+\frac{\theta}{256}\right)\rho}(\xi)\setminus B_{\left(\frac{3}{4}-\frac{\theta}{256}\right)\rho}(\xi).$$
To prove this notice that due to \thref{Monotonie}
\begin{equation*}
\theta^2(\mu,x)\le\frac{1}{4\pi}\Will(\Sigma)\le2-\frac{\delta_0}{4\pi}\quad\text{for all }x\in\Rdrei.
\end{equation*}
Now assume that our claim is false, i.e. there exists a sequence $\rho_l\to0$ such that
$$\left(\Sigma\cap B_{\left(\frac{3}{4}+\frac{\theta}{256}\right)\rho_l}(\xi)\setminus B_{\left(\frac{3}{4}-\frac{\theta}{256}\right)\rho_l}(\xi)\right)\setminus\mathcal{A}(\rho_l)\neq\emptyset\quad\text{for all }l.$$
Since we already know that $\Sigma$ can locally be written as a $C^{1,\alpha}\cap W^{2,2}$-graph away from the bad point $\xi$ we get that $\Sigma\cap B_{\frac{3}{4}\rho_1}(\xi)$ contains two components $\Sigma_1$ and $\Sigma_2$ such that $\Sigma_1\cap\Sigma_2=\{\xi\}$. Since $\Sigma_i$ can locally be written as a $C^{1,\alpha}\cap W^{2,2}$-graph in $B_{\frac{3}{4}\rho_1}(\xi)\setminus\{\xi\}$, we get that $\theta^2(\Sigma_i,x)=1$ for all $x\neq\xi$, and by upper semicontinuity that $\theta^2(\Sigma_i,\xi)\ge1$. Therefore it follows that $\theta^2(\mu,\xi)\ge\theta^2(\Sigma_1,\xi)+\theta^2(\Sigma_2,\xi)\ge2,$ a contradiction and the claim follows.\\[0,3cm]
From this and $\Sigma_k\to\Sigma$ we get for $\rho<\rho_0$ and $k\in\varmathbb{N}$ sufficiently large that
$$\Sigma_k\cap\mathcal{A}(\rho)\cap B_{\left(\frac{3}{4}+\frac{\theta}{512}\right)\rho}(\xi)\setminus B_{\left(\frac{3}{4}-\frac{\theta}{512}\right)\rho}(\xi)=\Sigma_k\cap B_{\left(\frac{3}{4}+\frac{\theta}{512}\right)\rho}(\xi)\setminus B_{\left(\frac{3}{4}-\frac{\theta}{512}\right)\rho}(\xi).$$
Define the set
$$C_k=\left\{s\in\left(0,\theta\frac{\rho}{1024}\right)\,\Bigg|\,\partial B_{\frac{3}{4}\rho+s}(\xi)\cap L\cap\bigcup_m d_{k,m}=\emptyset\right\}.$$
The diameter estimates for the discs $d_{k,m}$ yield for $\varepsilon\le\varepsilon_0$ that $\Leins(C_k)\ge\theta\frac{\rho}{2048}.$ The selection principle in \cite{SL} yields that there exists a set $C\subset\left(0,\theta\frac{\rho}{1024}\right)$ with $\Leins(C)\ge\theta\frac{\rho}{2048}$ and such that every $s\in C$ lies in $C_k$ for infinitely many $k\in\varmathbb{N}$.\\[0,3cm]
Now define the set
$$D_k=\left\{s\in C\,\Bigg|\,\int_{\graph {u_k}_{|_{\partial B_{\frac{3}{4}\rho+s}(\xi)\cap L}}}|\rm{A}_k|^2\,d\Hzwei\le\frac{4096}{\theta\rho}\int_{\Sigma_k\cap\mathcal{A}(\rho)}|\rm{A}_k|^2\,d\Hzwei\right\}.$$
By a simple Fubini-type argument (as done before) it follows that $\Leins(D_k)\ge\theta\frac{\rho}{4096}$, and again by the selection principle there exists a $s\in\left(0,\theta\frac{\rho}{1024}\right)$ such that $s\in D_k$ for infinitely many $k\in\varmathbb{N}$. It follows that $u_k$ is defined on the circle $\partial B_{\frac{3}{4}\rho+s}(\xi)\cap L$ and that $\graph {u_k}_{|_{\partial B_{\frac{3}{4}\rho+s}(\xi)\cap L}}$ divides $\Sigma_k$ into two connected topological discs $\Sigma_1^k, \Sigma_2^k$, one of them, w.l.o.g. $\Sigma_1^k$, intersecting $B_{\frac{3}{4}\rho}(\xi)$.\\[0,3cm]
From the estimates for the function $u_k$ and the choice of $s$ we have
$$\graph {u_k}_{|_{\partial B_{\frac{3}{4}\rho+s}(\xi)\cap L}}\subset\mathcal{A}(\rho)\cap B_{\left(\frac{3}{4}+\frac{\theta}{512}\right)\rho}(\xi)\setminus B_{\left(\frac{3}{4}-\frac{\theta}{512}\right)\rho}(\xi).$$
From this inclusion and
$$\Sigma_k\cap\mathcal{A}(\rho)\cap B_{\left(\frac{3}{4}+\frac{\theta}{512}\right)\rho}(\xi)\setminus B_{\left(\frac{3}{4}-\frac{\theta}{512}\right)\rho}(\xi)=\Sigma_k\cap B_{\left(\frac{3}{4}+\frac{\theta}{512}\right)\rho}(\xi)\setminus B_{\left(\frac{3}{4}-\frac{\theta}{512}\right)\rho}(\xi)$$
we get that
$$\Sigma_1^k\subset B_{\left(\frac{3}{4}+\frac{\theta}{512}\right)\rho}(\xi),$$
and the Monotonicity formula yields
$$\Hzwei\left(\Sigma_1^k\right)\le c\rho^2.$$
According to \thref{extension} let $w_k\in C^\infty\left(B_{\frac{3}{4}\rho+s}(\xi)\cap L,L^\perp\right)$ be an extension of $u_k$ restricted to $\partial B_{\frac{3}{4}\rho+s}(\xi)\cap L$. In view of the estimates for $u_k$ and therefore for $w_k$ we get that
$$\graph w_k\subset B_{\left(\frac{3}{4}+\frac{\theta}{512}\right)\rho}(\xi).$$
Now we can define the surface $\tilde\Sigma_k$ by
$$\tilde\Sigma_k=\Sigma_k\setminus\Sigma_{1}^k\cup\graph w_k.$$
By construction we have that $\tilde\Sigma_k$ is an embedded and connected $C^{1,1}$-surface with $\genus\tilde\Sigma_k=0$, which surrounds an open set $\tilde\Omega_k\subset\Rdrei$.\\[0,3cm]
The problem is again that $\tilde\Sigma_k$ might not be a comparison surface. But we can do the same correction as done before in \thref{2ff-absch} to get for all $\rho\le\rho_0$
\begin{equation}
\liminf_{k\to\infty}\int_{\Sigma_k\cap B_{\left(\frac{3}{4}-\frac{\theta}{512}\right)\rho}(\xi)}|\rm{A}_k|^2\,d\Hzwei\le c\rho^\alpha
\end{equation}
Thus by definition of the bad points $\xi$ could not have been a bad point and therefore the set of bad points is empty. Thus we have shown that for every point $\xi\in\Sigma$ there exists a radius $\rho>0$, a 2-dim. plane $L$ and a function
\begin{equation}
u\in C^{1,\alpha}\left(B_\rho(\xi)\cap L\right)\cap W^{2,2}\left(B_\rho(\xi)\cap L\right)
\end{equation}
for some $\alpha\in\left(0,\frac{1}{2}\right)$, with
\begin{equation}\label{500}
\int_{B_\tau(x)\cap L}|\D^2 u|^2\le c\tau^{2\alpha}
\end{equation}
for all $x\in B_\rho(\xi)\cap L$ and all $\tau<\rho$ such that $B_\tau(x)\subset B_\rho(\xi)\cap L$, and such that
\begin{equation}\label{400}
\Sigma\cap B_\rho(\xi)=\graph u\cap B_\rho(\xi).
\end{equation}
By definition of $\Sigma$ and an approximation argument we have
$$\Will(\Sigma)\le\inf_{\tilde\Sigma\in\M}\Will(\tilde\Sigma)=\inf_{\tilde\Sigma\in C^1\cap W^{2,2}, I(\tilde\Sigma)=\sigma}\Will(\tilde\Sigma).$$
On the other hand we have that the first variation of the isoperimetric ratio of $\Sigma$ is not equal to 0 as shown in \thref{vectorfield}. Therefore there exists a Lagrange multiplier $\lambda\in\Reins$ such that for all $\phi\in C_c^\infty((-\varepsilon,\varepsilon)\times\Rdrei,\Rdrei)$ with $\phi(0,\cdot)=0$
\begin{equation}
\frac{d}{dt}\Big(\Will(\phi_t(\Sigma))-\lambda I(\phi_t(\Sigma))\Big)_{|_{t=0}}=0.
\end{equation}
Restricting to $\phi\in C_c^\infty((-\varepsilon,\varepsilon)\times B_\rho(\xi),\Rdrei)$ and using the graph representation (\ref{400}) this yields after some computation that $u$ is a weak solution of
\begin{eqnarray}\label{50}
\partial_k\partial_l\left(\A_{ijkl}(\D u)\,\,\partial_i\partial_j u\right)+\partial_i\Be_i(\D u,\D^2 u)=\lambda\Big(\partial_i\Ce_i(\D u)+\Ce_0\Big)
\end{eqnarray}
for some coefficients $\A_{ijkl}, \Be_i, \Ce_i$ and $\Ce_0$ that perfectly fits into the scheme of Lemma 3.2 in \cite{SL}. Since by (\ref{500}) $u$ fulfills the assumptions of this lemma, we get by a bootstrap argument that $u$ is actually smooth.\\[0,3cm]
Therefore we have finally shown that $\Sigma$ can locally be written as a smooth graph and we get that $\Hzwei(\Sigma)=\mu(\Rdrei)=1$ and $\Sigma=\partial\Omega$. As mentioned before $\Omega$ has the right volume and therefore $\Sigma$ has the right isoperimetric ratio, especially $\Sigma\in\M$. Finally (\ref{4'}) yields that $\Sigma$ is a minimizer of the Willmore energy in the set $\M$ and therefore the existence part of \thref{MainThm} is proved. \\[0,3cm]
Last but not least we have to show that the function $\beta$ is continuous and strictly decreasing. For that let $0<\sigma_0<1$. Choose according to the above $\Sigma_0\in\mathcal{M}_{\sigma_0}$ such that $\Will(\Sigma_0)=\beta(\sigma_0)$. As in section 2 the Willmore flow $\Sigma_t$ with initial data $\Sigma_0$ exists smoothly for all times and converges to a round sphere. By a result of Bryant in \cite{BR}, which states that the only Willmore spheres with Willmore energy smaller than $8\pi$ are round spheres, it follows that $\Will(\Sigma_t)$ is strictly decreasing in $t$. Therefore for every $\sigma\in(\sigma_0,1]$ there exists a surface $\Sigma\in\M$ with $\Will(\Sigma)<\Will(\Sigma_0)=\beta(\sigma_0)$, and therefore $\beta(\sigma)<\beta(\sigma_0)$. To prove the continuity notice that the first variation of the isoperimetric ratio of $\Sigma_0$ is not equal to 0. As in \thref{2ff-absch}, where we corrected the isoperimetric ratio by applying a suitable variation, we can change the isoperimetric ratio of $\Sigma_0$ a little bit, in fact make it a little larger, without changing the $L^2$-norm of the second fundamental form, i.e. by Gauss-Bonnet the Willmore energy, too much. Therefore we get a new surface $\Sigma\in\M$ for a slightly larger $\sigma$ such that $|\Will(\Sigma_0)-\Will(\Sigma)|$ is small. Finally we get from the monotonicity of $\beta$ proved in section 2 
$$|\beta(\sigma)-\beta(\sigma_0)|=\beta(\sigma)-\beta(\sigma_0)\le\Will(\Sigma)-\Will(\Sigma_0).$$
This shows that $\beta$ is continuous and therefore \thref{MainThm} is now completely proved. \hfill$\square$

\section{Convergence to a double sphere}
In this section we prove the convergence to a double sphere stated in the introduction in \thref{Thm2}. For that let $\sigma_k\in(0,1)$ such that $\sigma_k\to0$. Choose according to \thref{MainThm} surfaces $\Sigma_k\in\mathcal{M}_{\sigma_k}$ such that $\Will(\Sigma_k)=\beta(\sigma_k)\le8\pi$. After scaling and translation we may assume that $0\in\Sigma_k$ and $\Hzwei(\Sigma_k)=1$. As in section 3 it follows (after passing to a subsequence) that
$$\mu_k=\Hzwei\llcorner\Sigma_k\to\mu\quad\text{in }C_c^0(\Rdrei)',$$
where $\mu$ is an integral, rectifiable 2-varifold in $\Rdrei$ with compact support, $\theta(\mu,\cdot)\ge1$ $\mu$-a.e. and weak mean curvature vector $\vec{H}_\mu\in L^2(\mu)$, such that
\begin{equation*}
\Will(\mu)\le\liminf_{k\to\infty}\Will(\Sigma_k)=\liminf_{k\to\infty}\beta(\sigma_k)=8\pi.
\end{equation*}
The last equation follows from \thref{MainThm}. Moreover we get as in section 3 that
$$\Sigma_k\to\spt\mu\quad\text{in the Hausdorff distance sense}.$$
Define again the bad points $\B$ with respect to a given $\varepsilon>0$ as in (\ref{7}).\\[0,3cm]
As before there exist only finitely many bad points and for every $\xi_0\in\spt\mu\setminus\B$ there exists a $\rho_0=\rho_0(\xi_0,\varepsilon)>0$ such that
\begin{equation*}
\int_{\Sigma_k\cap B_{\rho_0}(\xi_0)}|\rm{A}_k|^2\,d\Hzwei\le2\varepsilon^2\quad\text{for infinitely many }k\in\varmathbb{N}.
\end{equation*}
Let $\xi_0\in\spt\mu\setminus\B$ and choose a sequence $\xi_k\in\Sigma_k$ such that $\xi_k\to\xi_0$. For $k$ sufficiently large we may apply the graphical decomposition lemma to $\Sigma_k$, $\xi_k$ and $\rho<\frac{\rho_0}{2}$ to get for $\varepsilon\le\varepsilon_0$ that there exist pairwise disjoint closed subsets $P_1^k,\ldots,P_{N_k}^k$ of $\Sigma_k$ such that
$$\Sigma_k\cap\overline{B_\frac{\rho}{2}(\xi_k)}=\left(\bigcup_{j=1}^{J_k}\graph u_j^k\cup\bigcup_{n=1}^{N_k} P_n^k\right)\cap\overline{B_\frac{\rho}{2}(\xi_k)},$$
where the sets $P_n^k$ are topological discs disjoint from $\graph u_j^k$, $u_j^k\in C^\infty\left(\overline\Omega_{k,j},L_{k,j}^\perp\right)$, $\Omega_{k,j}=\left(B_{\lambda_{k,j}}(\pi_{L_{k,j}}(\xi_k))\cap L_{k,j}\right)\backslash\bigcup_{m=1}^{M_{k,j}}d_{k,j,m}$ with $\lambda_{k,j}>\frac{\rho}{2}$, $L_{k,j}$ is a 2-dim. plane and the sets $d_{k,j,m}$ are pairwise disjoint closed discs in $L_{k,j}$, and such that we have the estimates
$$\sum_m\diam d_{k,j,m}+\sum_n\diam P_n^k\le c\varepsilon^\frac{1}{2}\rho\quad\text{and}\quad\frac{1}{\rho}\|u_j^k\|_{L^\infty(\Omega_{k,j})}+\|Du_j^k\|_{L^\infty(\Omega_{k,j})}\le c\varepsilon^{\frac{1}{6}}.$$
We claim that for $\theta\in(0,1)$ and all $k$ sufficiently large (depending on $\rho$, $\theta$)
\begin{equation}\label{d1}
\graph u_j^k\cap B_{\theta\frac{\rho}{2}}(\xi_k)\neq\emptyset\quad\text{for at least two }j\in\{1,\ldots,J_k\}.
\end{equation}
Suppose this is false. Notice that at least one graph has to intersect with $B_{\theta\frac{\rho}{2}}(\xi_k)$ since $\xi_k\in\Sigma_k$ and because of the diameter estimates for the $P_n^k$'s. After passing to a subsequence we may assume that 
$$\Sigma_k\cap\overline{B_{\theta\frac{\rho}{2}}(\xi_k)}=\left(\graph u_k\cup\bigcup_{n=1}^{N_k} P_n^k\right)\cap\overline{B_{\theta\frac{\rho}{2}}(\xi_k)}$$
and $B_{\theta\frac{\rho}{4}}(\xi_0)\subset B_{\theta\frac{\rho}{2}}(\xi_k)$ for all $k$. Let $\chiup_k=\chiup_{\Omega_k}$, where $\Omega_k$ is the open set surrounded by $\Sigma_k$. Since the isoperimetric ratio $I(\Sigma_k)\to0$, it follows that $\chiup_k\to0$ in $L^1$. Let $g\in C_c^1(B_{\theta\frac{\rho}{4}}(\xi_0),\Rdrei)$. We get that
\begin{equation}\label{d2}
\int_{\Sigma_k}\big<g,\nu_k\big>\,d\Hzwei=\int\chiup_k\diver g\to0,
\end{equation}
where $\nu_k$ is the outer normal to $\partial\Omega_k=\Sigma_k$. By assumption we have
$$\int_{\Sigma_k}\big<g,\nu_k\big>\,d\Hzwei=\int_{\graph u_k\cap B_{\theta\frac{\rho}{4}}(\xi_0)}\big<g,\nu_k\big>\,d\Hzwei+\sum_n\int_{P_n^k\cap B_{\theta\frac{\rho}{4}}(\xi_0)}\big<g,\nu_k\big>\,d\Hzwei.$$
The Monotonicity formula and the diameter estimates yield that the second term is bounded by $c\varepsilon\rho^2$. Choose $g=\pm\varphi e_3$, where $\varphi\in C_c^1(B_{\theta\frac{\rho}{4}}(\xi_0))$ such that $\varphi\ge\chiup_{B_{\theta\frac{\rho}{8}}(\xi_0)}$. We get (by choosing the right sign and after rotation)
$$\int_{\graph u_k\cap B_{\theta\frac{\rho}{4}}(\xi_0)}\big<g,\nu_k\big>\,d\Hzwei\ge\int_{\left(B_{\lambda_k}(\pi_{L_k}(\xi_k))\cap L_k\right)\backslash\bigcup_m d_{k,m}}\chiup_{B_{\theta\frac{\rho}{8}}(\xi_0)}(x+u_k(x)).$$
It follows from the diameter estimates for the sets $P_n^k$ and the bounds on $u_k$ that $\dist(\xi_k,L_k)\le c\varepsilon^\frac{1}{6}\rho$. Since $\xi_k\to\xi_0$ we get for $\varepsilon\le\varepsilon_0$ and $k$ sufficiently large that $\chiup_{B_{\theta\frac{\rho}{8}}(\xi_0)}(x+u_k(x))=1$ if $x\in\left(B_{\theta\frac{\rho}{16}}(\pi_{L_k}(\xi_k))\cap L_k\right)\backslash\bigcup_m d_{k,m}$. The diameter estimates for the discs $d_{k,m}$ finally yield $\int_{\Sigma_k}\big<g,\nu_k\big>\,d\Hzwei\ge c\rho^2-c\varepsilon\rho^2$. In view of (\ref{d2}) we arrive for $\varepsilon\le\varepsilon_0$ at a contradiction.\\[0,3cm]
Now let $\rho<\frac{\rho_0}{2}$ such that $\mu(\partial B_\frac{\rho}{2}(\xi_0))=0$ and therefore $\mu_k(B_\frac{\rho}{2}(\xi_0))\to\mu(B_\frac{\rho}{2}(\xi_0))$. Let $\delta$, $\theta\in(0,\frac{1}{2})$. For $k$ sufficiently large we may assume that $B_{(1-\delta)\frac{\rho}{2}}(\xi_k)\subset B_\frac{\rho}{2}(\xi_0)$ and by (\ref{d1})
$$\graph u_1^k\cap B_{\theta(1-\delta)\frac{\rho}{2}}(\xi_k)\neq\emptyset\quad\text{and}\quad\graph u_2^k\cap B_{\theta(1-\delta)\frac{\rho}{2}}(\xi_k)\neq\emptyset.$$
Let $x_j^k\in\graph u_j^k\cap B_{\theta(1-\delta)\frac{\rho}{2}}(\xi_k)$. In view of the diameter estimates for the sets $P_n^k$ we get
$$\mu_k(B_\frac{\rho}{2}(\xi_0))\ge\sum_{j=1}^2\int_{\left(B_{\lambda_{k,j}}(\pi_{L_{k,j}}(\xi_k))\cap L_{k,j}\right)\backslash\bigcup_m d_{k,j,m}}\chiup_{B_{(1-\theta)(1-\delta)\frac{\rho}{2}}(x_j^k)}(x+u_j^k(x))-c\varepsilon\rho^2.$$
Since $x_j^k\in\graph u_j^k\cap B_{\theta(1-\delta)\frac{\rho}{2}}(\xi_k)$ we have that $x_j^k=z_j^k+u_j^k(z_j^k)$ with $z_j^k\in L_{k,j}$ such that $|z_j^k-\pi_{L_{k,j}}(\xi_k)|\le\theta(1-\delta)\frac{\rho}{2}$. Therefore $B_{(1-\theta)(1-\delta)\frac{\rho}{2}}(z_j^k)\subset B_{\lambda_{k,j}}(\pi_{L_{k,j}}(\xi_k))$. Moreover it follows from the bounds for $u_j^k$ that $\chiup_{B_{(1-\theta)(1-\delta)\frac{\rho}{2}}(x_j^k)}(x+u_j^k(x))=1$ if $|x-z_j^k|<\frac{(1-\theta)(1-\delta)}{1+c\varepsilon^\frac{1}{6}}\frac{\rho}{2}$. Therefore after all we get in view of the diameter estimates for the discs $d_{k,m,j}$ that
$$\mu_k(B_\frac{\rho}{2}(\xi_0))\ge2\left(\frac{(1-\theta)(1-\delta)}{1+c\varepsilon^\frac{1}{6}}\right)^2\pi\left(\frac{\rho}{2}\right)^2-c\varepsilon\rho^2\ge\frac{3}{2}\pi\left(\frac{\rho}{2}\right)^2$$
for $\varepsilon\le\varepsilon_0$ and $\delta$, $\theta$ sufficiently small. Thus for all $\xi_0\in\spt\mu\setminus\B$
$$\mu(B_\frac{\rho}{2}(\xi_0))\ge\frac{3}{2}\pi\left(\frac{\rho}{2}\right)^2.$$
Now since the density exists everywhere by \thref{Monotonie}, since $\mu$ is integral and since $\mu(\B)=0$ (which follows from the Monotonicity formula) we have shown that $\theta^2(\mu,\cdot)\ge2$ $\mu$-a.e.. Since $\Will(\mu)\le8\pi$, the Monotonicity formula in \thref{Monotonie} yields $2\le\theta^2(\mu,\cdot)\le\frac{1}{4\pi}\Will(\mu)\le2$ $\mu$-a.e. and therefore
\begin{equation*}
\theta^2(\mu,\cdot)=2\quad\mu\text{-a.e.}\quad\text{and}\quad\Will(\mu)=8\pi.
\end{equation*}
Now define the new varifold
$$\tilde\mu=\frac{1}{2}\mu.$$
It follows that $\tilde\mu$ is a rectifiable 2-varifold in $\Rdrei$ with compact support $\spt\tilde\mu=\spt\mu$ and weak mean curvature vector $\vec{H}_{\tilde\mu}=\vec{H}_\mu\in L^2(\tilde\mu)$, such that $\theta^2(\tilde\mu,\cdot)=1$ $\tilde\mu$-a.e. and $\Will(\tilde\mu)=4\pi$. The next lemma yields that $\tilde\mu$ is a round sphere in the sense that $\tilde\mu=\Hzwei\llcorner\partial B_r(a)$ for some $r>0$ and $a\in\Rdrei$. Therefore $\mu$ is a double sphere as claimed and \thref{Thm2} is proved.

\begin{Le}
Let $\mu\neq0$ be a rectifiable 2-varifold in $\Rdrei$ with compact support and weak mean curvature vector $\vec{H}\in L^2(\mu)$ such that
\begin{eqnarray*}
\hspace{-3cm}(i) & & \theta^2(\mu,x)=1\quad\text{for }\mu\text{-a.e. }x\in\Rdrei, \\
\hspace{-3cm}(ii) & & \Will(\mu)=\frac{1}{4}\int|\vec{H}|^2\,d\mu\le4\pi.
\end{eqnarray*}
Then $\mu$ is a round sphere, namely $\mu=\Hzwei\llcorner\partial B_r(a)$ for some $r>0$ and $a\in\Rdrei$.
\end{Le}

\begin{Pro}
From \thref{Monotonie} it follows that the density exists everywhere and that $\theta^2(\mu,x)\ge1$ for all $x\in\spt\mu$. But then \thref{Monotonie} yields 
\begin{equation}
\quad\Will(\mu)=4\pi\quad\text{and}\quad\theta^2(\mu,x)=1\quad\text{for all }x\in\spt\mu.
\end{equation}
Since $\mu\neq0$ it follows from \thref{Monotonie} that there exists a $R>0$ such that $\spt\mu\setminus B_R(x)\neq\emptyset$ for all $x\in\Rdrei$. Let $x_0\in\spt\mu$. Since $\spt\mu$ is compact it follows from \thref{Monotonie} that $|\vec{H}(x)|=4\left|\frac{(x-x_0)^\perp}{|x-x_0|^2}\right|\le\frac{8}{R}$ for $\mu$-a.e. $x\in\spt\mu\setminus B_\frac{R}{2}(x_0)$. On the other hand by choosing $x_1\in\spt\mu\setminus B_R(x_0)$ it follows that $|\vec{H}(x)|\le\frac{8}{R}$ for $\mu$-a.e. $x\in\spt\mu\setminus B_\frac{R}{2}(x_1)$. Since $B_\frac{R}{2}(x_0)\cap B_\frac{R}{2}(x_1)=\emptyset$ it follows that $|\vec{H}(x)|\le\frac{8}{R}$ for $\mu$-a.e. $x\in\spt\mu$ and therefore $\vec{H}\in L^\infty(\mu)$. Using Allard's regularity theorem (Theorem 24.2 in \cite{SLGMT}) we see that $\spt\mu$ can locally be written as a $C^{1,\alpha}$-graph $u$ for some $\alpha\in(0,1)$. As in (\ref{28}) it follows that $u$ is a weak solution of
$$\sum_{i,j=1}^2\partial_j\left(\sqrt{\det g}\,\,g^{ij}\partial_i F\right)=\sqrt{\det g}\,\,\vec{H}\circ F$$
where $F(x)=x+u(x)$, $g_{ij}=\delta_{ij}+\partial_i u\cdot\partial_j u$. Since $\vec{H}\in L^p(\mu)$ for every $p\ge1$, it follows from a standard difference quotient argument (as for example in \cite{GT}, Theorem 8.8) that $u\in W^{2,p}$ for every $p\ge1$ and therefore
\begin{equation}\label{101}
\int_{B_\rho}|\D^2 u|^2\le c\rho^\alpha.
\end{equation}
From a classical result of Willmore \cite{W} and an approximation argument we get
$$\Will(\mu)=4\pi\le\inf_{\text{smooth }\Sigma}\Will(\Sigma)=\inf_{\Sigma\in C^1\cap W^{2,2}}\Will(\Sigma).$$
Therefore $u$ solves the Euler-Lagrange equation (\ref{50}) (but with $\lambda=0$ since we do not have any constraints) and again with the power decay in (\ref{101}) and Lemma 3.2 in \cite{SL} it follows that $u$ is smooth. Thus $\spt\mu$ is a smooth surface with Willmore energy $4\pi$ and therefore a round sphere due to Willmore \cite{W}.
\end{Pro}

\renewcommand\thesection{\Alph{section}}\setcounter{section}{0}
\vspace{0,2cm}
\textbf{\Large{Appendix}}
\vspace{-0,3cm}
\section{Monotonicity formula}
Following L. Simon \cite{SL} and Kuwert/Sch\"{a}tzle \cite{KS} we state here a Monotonicity formula for rectifiable 2-variolds $\mu$ in $\Rdrei$ with square integrable weak mean curvature vector $\vec{H}\in L^2(\mu)$. We use the notation
$$\theta^2_\ast(\mu,\infty)=\liminf_{\rho\to\infty}\frac{\mu(B_\rho(0))}{\pi\rho^2},\quad\Will(\mu,E)=\frac{1}{4}\int_E|\vec{H}|^2\,d\mu\quad\text{for }E\subset\Rdrei\text{ Borel}.$$

\begin{Thm}\label{Monotonie}
Assume that $\vec{H}(x)\perp T_x\mu$ for $\mu$-a.e. $x\in\Rdrei$. Then the density
$$\theta^2(\mu,x)=\lim_{\rho\to0}\frac{\mu(B_\rho(x))}{\pi\rho^2}\quad\text{exists for all }x\in\Rdrei$$
and the function $\theta^2(\mu,\cdot)$ is upper semicontinuous. Moreover if $\theta^2_\ast(\mu,\infty)=0$, then we have for all $x_0\in\Rdrei$ and all $0<\sigma<\rho$
\begin{equation*}
\mu(B_\rho(x_0))\le c\rho^2,\phantom{\left|\frac{1}{4}\vec{H}(x)\right|}
\end{equation*}
\begin{equation*}
\theta^2(\mu,x_0)\le c\left(\frac{\mu(B_\rho(x_0))}{\pi\rho^2}+\Will(\mu,B_\rho(x_0))\right),\phantom{\left|\frac{1}{4}\vec{H}(x)\right|}
\end{equation*}
\begin{equation*}
\int_{B_\rho(x_0)\setminus B_\sigma(x_0)}\left|\frac{1}{4}\vec{H}(x)+\frac{(x-x_0)^\perp}{|x-x_0|^2}\right|^2\,d\mu(x)\le\frac{1}{4\pi}\Will(\mu)-\theta^2(\mu,x_0),
\end{equation*}
where $\perp$ denotes the projection onto $T_x\mu$.$\phantom{\left|\frac{1}{4}\vec{H}(x)+\frac{(x-x_0)^\perp}{|x-x_0|^2}\right|^2}$
\end{Thm}

\begin{Rem}
\rm{Brakke proved in chapter 5 of \cite{B} that $\vec{H}$ is perpendicular for any integral varifold with locally bounded first variation. Therefore the statements of this section apply to integral varifolds with square integrable weak mean curvature vector.}
\end{Rem}

\section{The graphical decomposition lemma of L. Simon}
Here we state the graphical decomposition lemma of Simon proved in \cite{SL}.
\begin{Thm}\label{Decomposition}
Let $\Sigma\subset\varmathbb R^n$ be a smooth surface. For given $\xi\in\Sigma$ and $\rho>0$ let
\begin{eqnarray*}
 & (i) & \partial\Sigma\cap\overline{B_\rho(\xi)}=\emptyset,\phantom{\int_{\Sigma\cap B_\rho(0)}} \\
 & (ii) & \Hzwei\left(\Sigma\cap\overline{B_\rho(\xi)}\right)\le\beta\rho^2\quad\text{for some }\beta>0,\phantom{\int_{\Sigma\cap                B_\rho(0)}} \\
 & (iii) & \int_{\Sigma\cap B_\rho(\xi)}|\textnormal{A}|^2\,d\Hzwei\le\varepsilon^2.
\end{eqnarray*}
Then there exists a $\varepsilon_0=\varepsilon_0(n,\beta)>0$ such that if $\varepsilon\le\varepsilon_0$ there exist pairwise disjoint closed subsets $P_1,\ldots,P_N$ of $\Sigma$ such that
$$\Sigma\cap\overline{B_\frac{\rho}{2}(\xi)}=\left(\bigcup_{j=1}^J\graph u_j\cup\bigcup_{n=1}^N P_n\right)\cap\overline{B_\frac{\rho}{2}(\xi)},$$
where the following holds:
\begin{enumerate}
\item The sets $P_n$ are topological discs disjoint from $\graph u_j$.
\item $u_j\in C^\infty\left(\overline\Omega_j,L_j^\perp\right)$, where $L_j\subset\varmathbb R^n$ is a 2-dim. plane and
$$\Omega_j=\left(B_{\lambda_j}(\pi_{L_j}(\xi))\cap L_j\right)\backslash\bigcup_{m=1}^M d_{j,m},$$
where $\lambda_j>\frac{\rho}{2}$ and the sets $d_{j,m}$ are pairwise disjoint closed discs in $L_j$.
\item Let $\tau\in\left(\frac{\rho}{4},\frac{\rho}{2}\right)$ such that $\Sigma\cap\partial B_\tau(\xi)$ is transversal and $\partial B_\tau(\xi)\cap\left(\bigcup_{n=1}^N P_n\right)=\emptyset$. Denote by $\{\Sigma_l\}_{l=1}^L$ the components of $\Sigma\cap B_\frac{\rho}{2}(\xi)$ such that $\Sigma_l\cap\overline{B_{\frac{\rho}{8}}(\xi)}\neq\emptyset$. It follows (after renumeration) that
$$\Sigma_l\cap\overline{B_\tau(\xi)}=D_{\tau,l}=\left(\graph u_l\cup\bigcup_{n=1}^N P_n\right)\cap\overline{B_\tau(\xi)},$$
where $D_{\tau,l}$ is a topological disc.
\item The following inequalities hold:
\begin{eqnarray*}
& & \sum_{m=1}^M\diam d_{j,m}\le c(n)\left(\int_{\Sigma\cap B_\rho(\xi)}|\rm{A}|^2\,d\Hzwei\right)^\frac{1}{4}\rho\le                    c(n)\varepsilon^\frac{1}{2}\rho, \\
& & \sum_{n=1}^N\diam P_n\le c(n,\beta)\left(\int_{\Sigma\cap B_\rho(\xi)}|\rm{A}|^2\,d\Hzwei\right)^\frac{1}{4}\rho\le                    c(n,\beta)\varepsilon^\frac{1}{2}\rho, \\
& & \frac{1}{\rho}\|u_j\|_{L^\infty(\Omega_j)}+\|Du_j\|_{L^\infty(\Omega_j)}\le c(n)\varepsilon^{\frac{1}{2(2n-3)}}.\phantom{\left(\int_{\Sigma}|\rm{A}|^2\,d\Hzwei\right)^\frac{1}{4}}
\end{eqnarray*}
\end{enumerate}
\end{Thm}

\section{Useful results}
\rm{In this section we state some useful results we need for the proof of \thref{MainThm}. \thref{extension} is an extension result adapted to the cut-and-paste procedure we use and is proved in \cite{JS}.}
\begin{Le}\label{extension}
Let $L$ be a 2-dim. plane in $\varmathbb R^n$, $x_0\in L$ and $u\in\C^\infty\left(U,L^\perp\right)$, where $U\subset L$ is an open neighborhood of $L\cap\partial B_\rho(x_0)$. Moreover let $|\D u|\le c$ on $U$. Then there exists a function $w\in\C^\infty(\overline{B_\rho(x_0)},L^\perp)$ such that
\begin{eqnarray*}
& (i) & w=u\quad\text{and}\quad\frac{\partial w}{\partial\nu}=\frac{\partial u}{\partial\nu}\quad\text{on }\partial B_\rho(x_0),\phantom{\int_{B_\rho}} \\
& (ii) & \frac{1}{\rho}||w||_{L^\infty(B_\rho(x_0))}\le c(n)\left(\frac{1}{\rho}||u||_{L^\infty(\partial B_\rho(x_0))}+||\D u||_{L^\infty(\partial B_\rho(x_0))}\right),\phantom{\int_{B_\rho}} \\
& (iii) & ||\D w||_{L^\infty(B_\rho(x_0))}\le c(n)||\D u||_{L^\infty(\partial B_\rho(x_0))},\phantom{\int_{B_\rho}} \\
& (iv) & \int_{B_\rho(x_0)}|\D^2w|^2\le c(n)\rho\int_{\graph u_{|_{\partial B_\rho(x_0)}}}|\rm{A}|^2\,d\Heins.
\end{eqnarray*}
\end{Le}
\begin{Pro}
After translation and rotation we may assume that $x_0=0$ and $L=\Rzwei\times\{0\}$. Moreover we may assume that $\rho=1$, the general result follows by scaling. \\[0,3cm]
Let $\phi\in\C^\infty(\overline{B_1(0)})$ be a cutoff-function such that $0\le\phi\le1$, $\phi=1$ on $B_\frac{1}{2}(0)$, $\phi=0$ on $\overline{B_1(0)}\backslash B_\frac{3}{4}(0)$ and $|\D\phi|+|\D^2\phi|\le c(n)$, and define the function $w_1\in C^\infty(\overline{B_1(0)})$ by
$$w_1(x)=\left(1-\phi(x)\right)u\left(\frac{x}{|x|}\right)+\phi(x)\fint_{\partial B_1(0)}u.$$ 
It follows that
$$w_1=u,\quad\frac{\partial w_1}{\partial\nu}=0\quad\text{on }\partial B_1(0),$$
$$||w_1||_{L^\infty(B_1(0))}\le c(n)||u||_{L^\infty(\partial B_1(0))},\quad||\D w_1||_{L^\infty(B_1(0))}\le c(n)||\D u||_{L^\infty(\partial B_1(0))}.$$
Using the Poincar\'e-inequality we also get
$$\int_{B_1(0)}|\D^2w_1|^2\le c(n)||u||^2_{W^{2,2}(\partial B_1(0))}.$$
Next let $w_2\in\C^\infty(\overline{B_1(0)})$ be the unique solution of the elliptic boundary value problem given by
$$\Delta w_2=0\quad\text{in }B_1(0),\quad w_2=\frac{\partial u}{\partial\nu}\quad\text{on }\partial B_1(0).$$
The solution $w_2$ is explicitly given by
$$w_2(x)=\frac{1}{2\pi}\int_{\partial B_1(0)}\frac{1-|x|^2}{|x-y|^2}\frac{\partial u}{\partial\nu}(y)\,dy.$$
Using standard estimates it follows that 
$$||w_2||_{L^\infty( B_1(0))}\le||\D u||_{L^\infty(\partial B_1(0))},\quad|\D w_2(x)|\le\frac{6}{1-|x|^2}||\D u||_{L^\infty(\partial B_1(0))},$$
$$||w_2||^2_{W^{1,2}(B_1(0))}\le c(n)\left(||\D u||^2_{L^2(\partial B_1(0))}+||\D^2 u||^2_{L^2(\partial B_1(0))}\right).$$
Next let $w_3\in\C^\infty(\overline{B_1(0)})$ be given by
$$w_3(x)=\frac{1}{2}\left(|x|^2-1\right)w_2(x).$$
It follows that
$$w_3=0,\quad\frac{\partial w_3}{\partial\nu}(x)=w_2(x)=\frac{\partial u}{\partial\nu}(x)\quad\text{on }\partial B_1(0),\phantom{\frac{\partial w_3}{\partial\nu}}$$
$$||w_3||_{L^\infty(B_1(0))}\le c||w_2||_{L^\infty(B_1(0))}\le c||\D u||_{L^\infty(\partial B_1(0))},\phantom{\frac{\partial w_3}{\partial\nu}}$$
$$||\D w_3||_{L^\infty(B_1(0))}\le c||\D u||_{L^\infty(\partial B_1(0))}.\phantom{\frac{\partial w_3}{\partial\nu}}$$
Moreover
$$\Delta w_3(x)=w_2(x)+x\cdot\D w_2(x)\quad\text{in }B_1(0).$$
Using again standard estimates it follows that
$$\int_{B_1(0)}|\D^2w_3|^2\le c\left(||\D u||^2_{L^2(\partial B_1(0))}+||\D^2 u||^2_{L^2(\partial B_1(0))}\right).$$
Finally define $w\in\C^\infty(\overline{B_1(0)})$ by
$$w(x)=w_1(x)+w_3(x).$$
The properties of $w_1$ and $w_3$ yield
$$w=u,\quad\frac{\partial w}{\partial\nu}=\frac{\partial u}{\partial\nu}\quad\text{on }\partial B_1(0),\phantom{\frac{\partial w}{\partial\nu}}$$
$$||w||_{L^\infty(B_1(0))}\le c\left(||u||_{L^\infty(\partial B_1(0))}+||\D u||_{L^\infty(\partial B_1(0))}\right),\phantom{\frac{\partial w}{\partial\nu}}$$
$$||\D w||_{L^\infty(B_1(0))}\le c||\D u||_{L^\infty(\partial B_1(0))},\phantom{\frac{\partial w}{\partial\nu}}$$
$$\int_{B_1(0)}|\D^2w|^2\le c||u||^2_{W^{2,2}(\partial B_1(0))}.\phantom{\frac{\partial w}{\partial\nu}}$$
By subtracting an appropriate linear function from $w$, using again the Poincar\'e-inequality and the assumption $|\D u|\le c$ we can get a better estimate for the $L^2$-norm of $\D^2w$, namely
$$\int_{B_1(0)}|\D^2w|^2\le c\int_{\partial B_1(0)}|\D^2u|^2\le c\int_{\graph u_{|\partial B_1(0)}}|\rm{A}|^2,$$
and the lemma is proved.
\end{Pro}
The second lemma is a decay result we need to get a power decay for the $L^2$-norm of the second fundamental form.
\begin{Le}\label{decay}
Let $g:(0,b)\to[0,\infty)$ be a bounded function such that 
$$g\left(x\right)\le\gamma g(2x)+cx^\alpha\quad\text{for all }x\in\left(0,\frac{b}{2}\right),$$
where $\alpha>0$, $\gamma\in(0,1)$ and $c$ some positive constant. There exists a $\beta\in(0,1)$ and a constant $c=c\left(b,||g||_{L^\infty(0,b)}\right)$ such that
$$g(x)\le cx^\beta\quad\text{for all }x\in\left(0,b\right).$$
\end{Le}
\vspace{0,5cm}
The last statement is a generalized Poincar\'{e} inequality proved by Simon in \cite{SL}.

\begin{Le}\label{Poincare}
Let $\mu>0$, $\delta\in\left(0,\frac{\mu}{2}\right)$ and $\Omega=B_\mu(0)\backslash E$, where $E$ is measurable with $\Leins(p_1(E))\le\frac{\mu}{2}$ and $\Leins(p_2(E))\le\delta$ where $p_1$ is the projection onto the $x$-axis and $p_2$ is the projection onto the $y$-axis. Then for any $f\in C^1(\Omega)$ there exists a point $(x_0,y_0)\in\Omega$ such that 
$$\int_\Omega\left|f-f(x_0,y_0)\right|^2\le C\mu^2\int_\Omega\left|\D f\right|^2+C\delta\mu\sup_\Omega|f|^2,$$
where $C$ is an absolute constant.
\end{Le}

\end{document}